\documentclass[12pt]{article}
\usepackage{amssymb}
\usepackage{mathrsfs}
\usepackage[top=72pt,bottom=96pt,left=72pt,right=72pt]{geometry}
\def\1{\mathbf{1}}
\def\#{\sharp}
\def\A{\mathbb{A}}
\def\Ad{\mathrm{Ad}}

\def\Aut{\mathrm{Aut}}
\def\b{\flat}
\def\C{\mathbb{C}}
\def\Cas{\mathrm{Cas}}
\def\Cl{\mathrm{Cl}}

\def\End{\mathrm{End}}
\def\e{\varepsilon}
\def\E{\mathbf{E}}
\def\EH{\mathcal{S}}
\def\ev{\mathrm{ev}}
\def\ext{\mathrm{ext}}
\def\F{\mathbf{F}}
\def\g{\mathfrak{g}}
\def\Gr{\mathrm{Gr}}
\def\H{\mathbb{H}}
\def\h{\mathfrak{h}}
\def\Hom{\mathrm{Hom}}
\def\id{\mathrm{id}}
\def\Im{\mathrm{Im}}
\def\ins{\lrcorner}
\def\L{\Lambda}
\def\m{\mathfrak{m}}
\def\M{\mathscr{M}}
\def\N{\mathbb{N}}
\def\O{\mathbb{O}}

\def\proof{\pfill\textbf{Proof:}\quad}
\def\pfill{\vskip6pt plus3pt minus2pt\noindent}
\def\qed{\ensuremath{\hfill\Box}\pfill}
\def\R{\mathbb{R}}
\def\Re{\mathrm{Re}\,}
\def\Ric{\mathrm{Ric}}
\def\S{\mathrm{Sym}}
\def\scal{\mathrm{scal}}
\def\sl{\mathfrak{sl}}
\def\so{\mathfrak{so}}
\def\SO{\mathbf{SO}}
\def\sp{\mathfrak{sp}}
\def\Sp{\mathbf{Sp}}
\def\Spin{\mathbf{Spin}}
\def\SU{\mathbf{SU}}
\def\t{\mathfrak{t}}
\def\tr{\mathrm{tr}}
\def\U{\mathbf{U}}
\def\vol{\mathrm{vol}}
\def\Z{\Sigma}
\def\<#1,#2>{\langle\,#1,\,#2\,\rangle}
\newtheorem{Lemma}{Lemma}[section]
\newtheorem{Remark}[Lemma]{Remark}
\newtheorem{Theorem}[Lemma]{Theorem}

\newtheorem{Corollary}[Lemma]{Corollary}
\newtheorem{Definition}[Lemma]{Definition}

\date{December 13th, 2020}
\title{Stability of Compact Symmetric Spaces}
\author{Uwe Semmelmann\footnote{Institut f\"ur Geometrie und Topologie,
 Fachbereich Mathematik, Universit\"at Stuttgart, Pfaffenwaldring 57,
 70569 Stuttgart, ALLEMAGNE.}\ \ 
 \& Gregor Weingart\footnote{Instituto de Matem\'aticas, Universidad
 Nacional Aut\'onoma de M\'exico, Avenida Universidad s/n, Colonia
 Lomas de Chamilpa, 62210 Cuernavaca, MEXIQUE.}}
\begin{document}
\maketitle
\begin{abstract}
 \noindent
 In this article we study the stability problem for the Einstein--Hilbert
 functional on compact symmetric spaces following and completing the seminal
 work of Koiso on the subject. We classify in detail the irreducible
 representations of simple Lie algebras with Casimir eigenvalue less than
 the Casimir eigenvalue of the adjoint representation, and use this
 information to prove the stability of the Einstein metrics on both
 the quaternionic and Cayley projective plane. Moreover we prove that the
 Einstein metrics on quaternionic Grassmannians different from projective
 spaces are unstable.
 \\[5pt]
 \begin{tabular}{lp{250pt}}
  \textit{MSC (2000):}& 53C25, 53C27, 53C44.
  \\
  \textit{Keywords:}&  Symmetric spaces, Einstein metrics, stability.
 \end{tabular}
\end{abstract}
\section{Introduction}
\label{intro}
 Perhaps the most interesting mathematical insight gained from studying general
 relativity is that Einstein metrics $g$ on a closed manifold $M$ can be
 characterized variationally as critical points of the Einstein--Hilbert
 or total scalar curvature functional
 $$
  \EH[\,g\,]
  \;\;:=\;\;
  \int_M\scal_g\,|\,\vol_g\,|\ .
 $$
 under volume preserving variations. In order to define the Hessian $\EH''$
 of the Einstein--Hilbert functional $\EH$ in an Einstein metric $g$ we
 recall that the set $\M\,:=\,\Gamma(\,\S^2_{\mathrm{reg}}T^*M\,)$ of
 pseudo--Riemannian metrics on $M$ is a Fréchet symmetric space under
 the binary operation
 $$
  *:\;\;\M\;\times\;\M\;\longrightarrow\;\M,\qquad
  (\,g,\,h\,)\;\longmapsto\;g\,*\,h\ ,
 $$
 characterized by $\#_{g\,*\,h}\,:=\,\#_g\,\circ\,\b_h\,\circ\,\#_g$,
 where $\b_g:\,TM\longrightarrow T^*M,\,X\longmapsto g(\,X,\,\cdot\,),$
 and $\#_g\,:=\,\b^{-1}_g$ are the musical isomorphisms associated to $g$.
 Geodesics in $\M$ are by definition continuous maps $\R\longrightarrow\M,
 \,t\longmapsto g_t,$ satisfying $g_s\,*\,g_t\,=\,g_{2s-t}$ for all
 $s,\,t\,\in\,\R$. The binary operation $*$ on $\M$ is evidently defined
 pointwise, in turn all geodesics can be written in the form $g_t\,=\,g_0
 (\,e^{tF}\,\cdot\,,\;\cdot\,)$ with a symmetric endomorphism field
 $F\,\in\,\Gamma(\,\End\,TM\,)$ with respect to $g_0\,\in\,\M$. In
 particular the exponential map in an Einstein metric $g\,\in\,\M_+$
 \begin{equation}\label{gexp}
  \exp_g:\;\;\Gamma(\,\S^2T^*M\,)\;\longrightarrow\;\M_+,\qquad
  h\;\longmapsto\;g(\;e^{\#_g\,\circ\,\b_h}\;\cdot\,,\;\cdot\;)\ ,
 \end{equation}
 is a diffeomorphism to the Fréchet manifold $\M_+\,\subseteq\,\M$
 of Riemannian metrics on $M$. Unluckily the level sets of the volume
 functional on $\M_+$ are not totally geodesic submanifolds of $\M_+$,
 nevertheless every geodesic through $g\,=\,g_0$ can be written as
 conformal change
 $$
  \exp_g(\;th\;)
  \;\;=\;\;
  e^{\frac tm\,\tr_gh}
  \;\exp_g(\;th_\circ\;)
  \qquad\qquad
  h_\circ
  \;\;:=\;\;
  h\;-\;\frac1m\,(\,\tr_gh\,)\,g
 $$
 of a volume preserving geodesic with $m\,:=\,\dim\,M$, $\tr_gh\,:=\,
 \tr(\,\#_g\circ\b_h\,)\,\in\,C^\infty(\,M\,)$ and a section $h_\circ\,\in\,
 \Gamma(\,\S^2_\circ T^*M\,)$ of the subbundle $\S^2_\circ T^*M\,\subseteq\,
 \S^2T^*M$ of trace free symmetric $2$--tensors with respect to $g$. Moreover
 the diffeomorphism group $\mathrm{Diff}\,M$ acts via pull back on the
 symmetric space $\M_+$ and preserves the volume functional on $\M_+$.
 In order to address both problems we decompose the domain of $\exp_g$
 for every compact Einstein manifold $(\,M,\,g\,)$ different from a round
 sphere into the direct sum (cf.~\cite{B87}, Lemma 4.57)
 \begin{equation}\label{deco}
  \Gamma(\,\S^2T^*M\,)
  \;\;=\;\;
  C^\infty(\,M\,)\,g\;\oplus\;\Im\;D\;\oplus\;
  \bigg(\;\ker\,D^*\,\cap\,\Gamma(\,\S^2_\circ T^*M\,)\;\bigg)\ ,
 \end{equation}
 where $D^*:\,\Gamma(\,\S^2T^*M\,)\longrightarrow\Gamma(\,TM\,),\,h\longmapsto
 D^*h,$ is the differential operator formally adjoint to the Lie derivative
 $D:\,\Gamma(\,TM\,)\longrightarrow\Gamma(\,\S^2T^*M\,),\,X\longrightarrow
 \mathfrak{Lie}_Xg$. Both differential operators $D$ and $D^*$ can be written
 in terms of the covariant derivatives $\nabla$ on the vector bundles
 $\S^2T^*M$ and $TM$ induced by the Levi--Civita connection on $M$, namely
 \begin{equation}\label{ll*}
  D X
  \;\;:=\;\;
  \sum_\mu\,E_\mu^\b\,\cdot\,(\,\nabla_{E_\mu}X\,)^\b
  \qquad\qquad
  D^*h
  \;\;:=\;\;
  -\;\sum_\mu\,(\,E_\mu\;\ins\;\nabla_{E_\mu}h\,)^\#\ ,
 \end{equation}
 where $E_1,\,\ldots,\,E_m\,\in\,\Gamma_{\mathrm{loc}}(\,TM\,)$ denotes
 some local orthonormal basis for $TM$. Up to sign $D^*h$ agrees with
 the divergence of the symmetric $2$--tensor $h\,\in\,\Gamma(\,\S^2T^*M\,)$.
 
 \pfill
 The stability problem for a given Einstein metric $g$ on a compact
 manifold $M$ different from the round metric on $S^m$ refers to the
 behaviour of the Einstein--Hilbert functional $\EH$ restricted to
 the image under $\exp_g$ of the third summand in the decomposition
 (\ref{deco}), the space $\ker\,D^*\,\cap\,\Gamma(\,\S^2_\circ T^*M\,)$
 of trace and divergence free symmetric $2$--tensors or tt--tensors on
 the Einstein manifold $M$. The Hessian of the Einstein--Hilbert functional
 in its critical point $g$ describes the second variation $\EH''$ of $\EH$
 along the geodesics in $\M_+$ through $g$ in the direction of tt--tensors
 $h_\circ\,\in\,\ker\,D^*\,\cap\,\Gamma(\,\S^2_\circ T^*M\,)$. It is given
 by (cf. \cite{B87}, Theorem 4.60)
 \begin{equation}\label{2var}
  \EH''[\,g;\;h_\circ\,]
  \;\;=\;\;
  -\;\frac12\;\int_Mg_{T^*\otimes T^*}\left(\;\Big[\;\Delta_L\;-\;
  2\,\frac\scal m\;\Big]\,h_\circ,\;h_\circ\;\right)\;|\,\vol_g\,|\ ,
 \end{equation}
 where the Lichnerowicz Laplacian $\Delta_L$ is a Laplace type operator
 defined by: 
 \begin{equation}\label{lichl}
  \Delta_L
  \;\;:=\;\;
  \nabla^*\nabla\;+\;2\,q(R)
  \qquad\qquad
  q(\,R\,)
  \;\;:=\;\;
  \frac14\sum_{\mu\nu}(E_\mu\wedge E_\nu)\,\star\, R_{E_\mu,\,E_\nu}\ .
 \end{equation}
 The interested reader may find more details on $\Delta_L$ and $q(\,R\,)$
 in \cite{SW19}, but should take notice that the curvature term $q(\,R\,)$
 is subject to a different normalization in that reference. The original
 definition of Lichnerowicz spelt out the curvature term $q(R)$ in the form
 $$
  2\,q(\,R\,)
  \;\;=\;\;
  2\,\raise10pt\hbox to0pt{$\;\scriptscriptstyle\circ$\hss}R\;+\;\Ric\ ,
 $$
 where $\Ric$ acts on symmetric tensors by $(\,\Ric\,h\,)(X,Y)\,:=\,
 h(\Ric\,X,Y)\,+\,h(X,\Ric\,Y)$ and:
 $$
  (\,\raise10pt\hbox to0pt{$\;\scriptscriptstyle\circ$\hss}R\,h\,)(\;X,\;Y\;)
  \;\;=\;\; 
  \sum_\mu h(\;R_{X,\,E_\mu}Y,\;E_\mu\;)\ .
 $$
 The main advantage of this way to define the curvature term $2\,q(\,R\,)$
 is that the integrand in equation (\ref{2var}) becomes $g_{T^*\otimes T^*}
 (\,[\,\nabla^*\nabla\,-\,2\raise10pt\hbox to0pt
 {$\;\scriptscriptstyle\circ$\hss}R\,]\,h_\circ,\,h_\circ\,)$.

 With $\Delta_L$ being a differential operator of Laplace type the Hessian
 $\EH''$ in equation (\ref{2var}) is non--positive except for possibly a
 finite dimensional subspace $R\,\subseteq\,\ker\,D^*\,\cap\,\Gamma(\,
 \S^2_\circ T^*M\,)$ of tt--tensors. Following the seminal work
 of Koiso we call an Einstein metric $g$ stable with respect to the
 Einstein--Hilbert functional $\EH$, if its Hessian $\EH''$ is negative
 on the space of tt--tensors; in the same vein we call $g$ unstable
 provided there are tt--tensor directions on which $\EH''$ is positive
 (cf.~\cite{Koiso80}, Definition 2.7). The kernel of $\EH''$ on tt--tensors
 agrees with the space of infinitesimal Einstein deformations characterized
 by the linearized Einstein equation
 $$
  \Delta_Lh_\circ\;-\;2\,\frac\scal m\,h_\circ
  \;\;=\;\;
  0\ ,
 $$
 where $\scal$ is the scalar curvature and $m\,:=\,\dim\,M$. In passing
 we observe that for Einstein metrics the constant $2\,\frac\scal m$ is
 exactly the eigenvalue of the Casimir operator on Killing vector fields.
 Extending the definitions of $\Delta_L$, $D$ and $D^*$ as well as the
 curvature term $q(\,R\,)$ to symmetric tensors $h\,\in\,\Gamma(\,\S^rT^*M\,)$
 of arbitrary degree $r\,\in\,\N$ we find in general $\Delta_Lh\,\geq\,
 4\,q(R)\,h$ for all divergence free symmetric tensors $D^*h\,=\,0$ with
 equality exactly for so--called Killing tensors satisfying $D\,h\,=\,0$
 (see \cite{HMS16}, Section 6).
 
 Stability of Einstein metrics with respect to the Einstein--Hilbert
 functional $\EH$ has been extensively studied for example in \cite{Koiso80},%
 \cite{Koiso82},\cite{DWW05} and \cite{K17}. In \cite{Koiso80},\cite{Koiso82}
 Koiso essentially classified the stable symmetric spaces of compact type.
 More precisely Koiso's result stipulates that the Einstein metrics on simply
 connected irreducible symmetric spaces of compact type are stable unless the
 spaces belong to one of the following three categories with
 $n\,\geq\,3;\,r,\,s\,\geq\,2$:
 \begin{description}
  \item[Irreducible Symmetric Spaces with Infinitesimal Deformations:]
   $$
    \hbox to70pt{$\displaystyle\SU(\,n\,)$\hfill}
    \hbox to110pt{$\displaystyle\SU(\,n\,)/_{\displaystyle\SO(\,n\,)}$\hfill}
    \hbox to115pt{$\displaystyle\SU(\,2n\,)/_{\displaystyle\Sp(\,n\,)}$\hfill}
    \hbox to70pt{$\displaystyle\Gr_r\C^{r+s}$\hfill}
    \hbox to70pt{$\displaystyle\E_6/_{\displaystyle\F_4}$\hfill}
   $$
  \item[Unstable Irreducible Symmetric Spaces:]
   $$
    \hbox to70pt{$\displaystyle\Sp(\,r\,)$\hfill}
    \hbox to110pt{$\displaystyle\Sp(\,n\,)/_{\displaystyle\U(\,n\,)}$\hfill}
    \hbox to115pt{$\displaystyle\Gr^{\mathrm{or}}_2\R^5$\hfill}
    \hbox to70pt{\hfill}\hbox to70pt{\hfill}
   $$
  \item[Irreducible Symmetric Spaces with undecided Stability Status:]
   $$
    \hbox to70pt{$\displaystyle\Gr_r\H^{r+s}$\hfill}
    \hbox to110pt{$\displaystyle\H P^2$\hfill}
    \hbox to115pt{$\displaystyle\O P^2$\hfill}
    \hbox to70pt{\hfill}\hbox to70pt{\hfill}
   $$
 \end{description}
 In the work of Koiso the complex quadric $\Gr^{\mathrm{or}}_2\R^5\,=\,
 \SO(5)/\SO(3)\times\SO(2)$ in $\C P^4$ was actually overlooked as remarked
 by Cao \& He \cite{CH15}, however it is unstable according to results by
 Gasqui \& Goldschmidt \cite{GG96}. Stability status remained undecided by
 Koiso only for the quaternionic and Cayley projective plane $\H P^2$ and
 $\O P^2\,=\,\F_4/\Spin(9)$, and the quaternionic Grassmannians
 $\Gr_r\H^{r+s}\,=\,\Sp(r+s)/\Sp(r)\times\Sp(s)$ with $r,\,s\,\geq\,2$.
 Apparently this question has not been settled since, in the recent work
 of Cao \& He \cite[Table 2]{CH15} for example the stability status of
 these spaces are listed as unknown. Our main result fills this gap
 and clarifies the stability status for the remaining symmetric spaces
 of compact type:
 
 \begin{Theorem}[Stability of Quaternionic and Cayley Projective Plane]
 \hfill\label{mt}\break
  The Cayley projective plane $\O P^2\,=\,\F_4/\Spin(9)$ is stable in the
  sense of Koiso. The quaternionic Grassmannians $\Gr_r\H^{r+s}\,=\,\Sp(r+s)
  /\Sp(r)\times\Sp(s)$ of quaternionic subspaces of dimension $r$ in $\H^{r+s}$
  are unstable in the sense of Koiso for all parameters $r,\,s\,\geq\,2$,
  for $r\,=\,1$ or $s\,=\,1$ however they are stable. In particular
  $\H P^2\,=\,\Gr_1\H^3$ is stable.
 \end{Theorem}

 \noindent
 It should be pointed out that Koiso identified the Lichnerowicz
 Laplacian on symmetric spaces with a suitably normalized Casimir operator
 $\Cas$ and used this information to compute the first eigenvalue of $\Delta_L$
 on the space $\Gamma(\,\S^2_\circ T^*M\,)$ of trace free symmetric
 $2$--tensors. The first eigenvalues $\frac{2(r+s)}{(r+s+1)}\,\frac\scal m$
 and $\frac43\,\frac\scal m$ Koiso obtained for the quaternionic Grassmannians
 $\Gr_r\H^{r+s}$ and the Cayley projective plane $\O P^2$ are both below the
 critical value $2\,\frac\scal m$. Our contribution to the classification of
 Koiso solves the question, whether the corresponding trace free eigentensors
 $h_\circ\,\in\,\Gamma(\,\S^2_\circ T^*M\,)$ of the Lichnerowicz Laplacian
 $\Delta_L$ can be chosen to be divergence free $D^*h_\circ\,=\,0$ as well
 or not.

 Besides Koiso's notion there actually exists a weaker notion of stability
 of Einstein metrics, the so--called $\EH$--linear stability (cf.~\cite{W17},
 \cite{SWW20}), which allows for the presence of infinitesimal deformations.
 More precisely an Einstein metric is called $\EH$--linearly stable, if
 $\EH''$ is non--positive on the space of tt--tensors. According to the
 classification of Koiso the only two symmetric spaces of compact type with
 infinitesimal deformations and subcritical eigenvalues on
 $\Gamma(\,\S^2_\circ T^*M\,)$ are $\SU(\,n\,)$ and $\E_6/\F_4$. In a
 forthcoming paper based on our approach Schwahn \cite{S20} shows that
 both symmetric spaces are $\EH$--linearly stable, because their
 subcritical eigenvalues are not realizable by tt--tensors.
 
 \pfill
 In Section \ref{prodiff} the left regular representation on sections of
 homogeneous vector bundles over and use the Frobenius reciprocity to
 associate a family of linear maps, the prototypical differential operators,
 to an equivariant differential operator. The prototypical divergence
 operators associated to $D^*$ allow us to translate the stability
 problem into a problem in finite dimensional linear algebra. In Section
 \ref{critreps} we provide the details of the classification of the critical
 representations by Koiso in order to identify the subcritical eigenspaces
 for the quaternionic and Cayley projective plane and the quaternionic
 Grassmannians. In Section \ref{geocay} and \ref{eingr} we decide for
 the Cayley projective plane $\O P^2$ and the quaternionic Grassmannians
 respectively, whether the subcritical eigenvalues can be realized by
 tt--tensors or not.
\section{Prototypical Differential Operators}
\label{prodiff}
 Analysis on homogeneous spaces or Harmonic Analysis is a subtopic of
 differential geometry of particular elegance, because many of its problems
 can be translated into equivalent problems of linear algebra by means of an
 extensive dictionary of rules and prescriptions. In this section we focus
 on a particular concept in this dictionary, the prototypical differential
 operators, in order to formulate the linear algebra equivalent of the
 stability problem for compact symmetric spaces in Corollary \ref{lsrs}.
 Of particular importance in Sections \ref{geocay} and \ref{eingr} is the
 formula (\ref{proto}) for the prototypical divergence operator $D^*_R$ we
 discuss at the end of this section in a form geared to be easily accesible
 and self--contained.

 \pfill
 Let us recall that a homogeneous space is a manifold $M$ endowed with a
 transitive smooth action $G\times M\longrightarrow M,\,(\,g,\,p\,)
 \longmapsto g\,\centerdot\,p,$ of a Lie group $G$ with the associated group
 homomorphism $\mu:\,G\longrightarrow\mathrm{Diff}\;M,\,g\longmapsto\mu_g,$
 of left multiplication $\mu_g(\,p\,)\,:=\,g\centerdot p$. A homogeneous
 vector bundle over a homogeneous space is a vector bundle $VM$ over the
 manifold $M$ endowed with a smooth action $\star:\,G\times VM\longrightarrow
 VM,\,(\,g,\,v\,)\longmapsto g\,\star\,v,$ on its total space, which covers
 the action of $G$ on $M$ and is linear on fibers. Every homogeneous vector
 bundle $VM$ over $M$ gives rise to the infinite dimensional left regular
 representation
 $$
  L:\;\;G\;\times\;\Gamma(\,VM\,)\;\longrightarrow\;\Gamma(\,VM\,),
  \qquad(\,g,\,v\,)\;\longmapsto\;L_gv\ ,
 $$
 of the group $G$ on the vector space $\Gamma(\,VM\,)$ by means of
 $(\,L_gv\,)(\,p\,)\,:=\,g\,\star\,v(\,g^{-1}\centerdot p\,)$. Every tensor
 bundle or more generally every natural vector bundle $VM$ on a homogeneous
 space $M$ is automatically a homogeneous vector bundle; the characteristic
 smooth action $\star:\,G\times VM\longrightarrow VM$ on its total space $VM$
 is implicitly defined by stipulating the identity
 \begin{equation}\label{req}
  L_gv
  \;\;\stackrel!=\;\;
  \mu_{g^{-1}}^*v
 \end{equation}
 for all $g\,\in\,G$. The rather unexpected change $g\rightsquigarrow g^{-1}$
 in this identity is mandated by the contravariance $(\varphi\,\circ\,\psi)^*
 \,=\,\psi^*\,\circ\,\varphi^*$ of the pull back of sections of natural vector
 bundles. A left invariant differential operator on a homogeneous space $M$
 is a differential operator
 \begin{equation}\label{ldo}
  D:\;\;\Gamma(\,VM\,)\;\longrightarrow\;\Gamma(\,WM\,),\qquad
  v\;\longmapsto\;D\,v\ ,
 \end{equation}
 between the sections of homogeneous vector bundles $VM$ and $WM$ over $M$,
 which is in addition equivariant $D(\,L_gv\,)\,=\,L_g(\,Dv\,)$ under $G$ for
 all $v\,\in\,\Gamma(\,VM\,)$ and $g\,\in\,G$. 

 The fiber of a homogeneous vector bundle $VM$ in a chosen base point
 $p\,\in\,M$ is naturally a representation $V_pM$ of the stabilizer or
 isotropy subgroup $H\,:=\,\{\;h\,\in\,G\;|\;h\,\centerdot\,p\,=\,p\;\}
 \,\subseteq\,G$ by restricting the smooth action $\star$ to the
 submanifold $H\times V_pM\,\subseteq\,G\times VM$. In terms of this
 representation the evaluation map $\ev_p:\,\Gamma(\,VM\,)\longrightarrow
 V_pM,\,v\longmapsto v(\,p\,)$ is equivariant under the stabilizer subgroup
 $H\,\subseteq\,G$ due to the trivial identity $(\,L_hv\,)(\,p\,)\,=\,
 h\,\star\,v(\,p\,)$ for all $h\,\in\,H$ and $v\,\in\,\Gamma(\,VM\,)$.
 Postcomposition with $\ev_p$ thus induces a linear map
 \begin{equation}\label{frob}
  \Hom_G(\;R,\;\Gamma(\,VM\,)\;)\;\stackrel\cong\longrightarrow\;
  \Hom_H(\;R,\;V_pM\;),\qquad F\;\longmapsto\;\ev_p\,\circ\,F
 \end{equation}
 for every finite dimensional representation $(\,R,\,\star_R\,)$ of the group
 $G$. The well--known Frobenius reciprocity asserts that this linear map is
 a vector space isomorphism, whose inverse
 \begin{equation}\label{invfrob}
  \Hom_H(\;R,\;V_pM\;)\;\stackrel\cong\longrightarrow\;
  \Hom_G(\;R,\;\Gamma(\,VM\,)\;),\qquad F\;\longmapsto\;F^\ext
 \end{equation}
 reads $(\,F^\ext r\,)(\,g\centerdot p\,)\,=\,g\,\star\,F(\,g^{-1}\star_Rr\,)$
 for all $r\,\in\,R$ and $g\,\in\,G$. Using Frobenius reciprocity we can break
 up a left invariant differential operator $D:\,\Gamma(\,VM\,)\longrightarrow
 \Gamma(\,WM\,)$ into more manageable pieces, the prototypical differential
 operators $D_R$ associated to $D$ and a finite dimensional representation
 $R$ of the group $G$ by means of the commutative diagram
 \begin{equation}\label{cd}
  \begin{picture}(270,65)(0,0)
   \put(  0, 0){$\Hom_H(\;R,\;V_pM\;)$}
   \put(  0,50){$\Hom_G(\,R,\,\Gamma(VM)\,)$}
   \put(170, 0){$\Hom_H(\;R,\;W_pM\;)$\ ,}
   \put(170,50){$\Hom_G(\,R,\,\Gamma(WM)\,)$}
   \put(104, 4){\vector(+1, 0){61}}
   \put(124, 7){$\scriptstyle D_R$}
   \put(104,54){\vector(+1, 0){61}}
   \put(126,57){$\scriptstyle D\,\circ$}
   \put( 54,45){\vector( 0,-1){34}}
   \put( 33,27){$\scriptstyle\ev_p\,\circ$}
   \put(224,45){\vector( 0,-1){34}}
   \put(230,27){$\scriptstyle\ev_p\,\circ$}
   \put(128,27){$\circlearrowleft$}
  \end{picture}
 \end{equation}
 where the upper arrow equals postcomposition with $D$ and the vertical
 arrows are the Frobenius reciprocity isomorphisms. Prototypical differential
 operators are used frequently in harmonic analysis, for example to calculate
 the spectra of elliptic left invariant differential operators on compact
 homogeneous spaces. In order to reduce both statements of our main Theorem
 \ref{mt} to a statement about prototypical differential operators let us call
 a finite dimensional $G$--invariant subspace $R\,\subseteq\,\Gamma(\,VM\,)$
 a $G$--characteristic subspace of $\Gamma(\,VM\,)$, if it contains the image
 $F(\,R\,)\,\subseteq\,R$ of every $G$--equivariant linear map
 $F:\,R\longrightarrow\Gamma(\,VM\,)$:

 \begin{Lemma}[Kernels of Left Invariant Differential Operators]
 \hfill\label{ked}\break
  Consider a left invariant differential operator $D:\,\Gamma(\,VM\,)
  \longrightarrow\Gamma(\,WM\,)$ between sections of homogeneous vector
  bundles $VM$ and $WM$ over a homogeneous space $M$ under a transitive
  smooth action $G\times M\longrightarrow M$ of a compact Lie group $G$.
  A finite dimensional $G$--characteristic subspace $R\,\subseteq\,
  \Gamma(\,VM\,)$ intersects the kernel of $D$ trivially $R\,\cap\,\ker\,D
  \,=\,\{\,0\,\}$, if and only if the prototypical differential operator
  associated to $D$ and $R$ is injective:
  $$
   D_R:\;\;\Hom_H(\;R,\;V_pM\;)\;\longrightarrow\;\Hom_H(\;R,\;W_pM\;)
  $$
 \end{Lemma}

 \proof
 Assuming for the moment that the prototypical differential operator $D_R$
 associated to $D$ and $R$ is injective, we observe that the kernel of the
 left invariant differential operator $D:\,\Gamma(\,VM\,)\longrightarrow
 \Gamma(\,WM\,)$ is necessarily a $G$--invariant subspace $\ker\,D\,\subseteq
 \,\Gamma(\,VM\,)$. The intersection $R\,\cap\,\ker\,D\,\subseteq\,R$ is thus
 a $G$--invariant subspace of the finite dimensional representation $R$ of
 the compact group $G$ so that there exists a surjective $G$--equivariant
 projection $P:\,R\longrightarrow R\,\cap\,\ker\,D$. Interpreting $P$ as a
 $G$--equivariant linear map $P:\,R\longrightarrow\Gamma(\,VM\,)$ we obtain
 an element $P$ of the kernel of postcomposition with the differential
 operator $D$:
 \begin{equation}\label{post}
  \Hom_G(\;R,\;\Gamma(\,VM\,)\;)\;\longrightarrow\;\Hom_G(\;R,\;
  \Gamma(\,WM\,)\;),\qquad F\;\longmapsto\;D\,\circ\,F\ .
 \end{equation}
 Postcomposition with $D$ however is conjugated to the injective prototypical
 differential operator $D_R$ via the commutative diagram (\ref{cd}) and so
 we conclude $P\,=\,0$ and in turn $R\,\cap\,\ker\,D\,=\,\{\,0\,\}$ due to
 the surjectivity of $P$. Conversely let us assume that $R$ and $\ker\,D$
 intersect trivially $R\,\cap\,\ker\,D\,=\,\{\,0\,\}$ and let $F:\,R
 \longrightarrow\Gamma(\,VM\,)$ be the $G$--equivariant linear map Frobenius
 reciprocal to a linear map $\ev_p\,\circ\,F\,\in\,\ker\,D_R$ in the kernel
 of the prototypical differential operator $D_R$. With $D\,\circ\,F\,=\,0$
 vanishing the image of $F$ is contained in $\ker\,D$, moreover $F(\,R\,)\,
 \subseteq\,R$ for the $G$--characteristic subspace $R$. In consequence
 $F(\,R\,)\,\subseteq\,R\,\cap\,\ker\,D\,=\,\{\,0\,\}$ leading to $F\,=\,0$
 and in turn to $\ker\,D_R\,=\,\{\,0\,\}$.
 \qed

 \pfill
 In the context of Theorem \ref{mt} we are interested in compact connected
 Riemannian symmetric spaces $M$ and the $G$--invariant subspace
 $R\,\subseteq\,\Gamma(\,\S^2_\circ T^*M\,)$ given by the direct
 sum of all eigenspaces of the Lichnerowicz Laplacian (\ref{lichl})
 restricted to trace free symmetric $2$--tensors
 $$
  \Delta_L:\;\;
  \Gamma(\,\S^2_\circ T^*M\,)\;\longrightarrow\;\Gamma(\,\S^2_\circ T^*M\,)
 $$
 for the eigenvalues below the critical value $2\,\frac\scal m$, where
 $\scal$ and $m$ denote the scalar curvature and the dimension of $M$
 respectively. Elliptic regularity \cite{LM} enjoyed by the elliptic
 differential operators of Laplace type like $\Delta_L$ ensures that
 $R\,\subseteq\,\Gamma(\,\S^2_\circ T^*M\,)$ is a finite dimensional
 $G$--invariant subspace. In order to verify that $R$ is a $G$--characteristic
 subspace we observe that the Lichnerowicz Laplacian $\Delta_L$ is a special
 case of the standard Laplace operator $\Delta$ defined in \cite{SW19} and
 thus agrees on every symmetric space $M$ with the Casimir operator of $G$.
 More precisely let $\Cas\,\in\,\mathscr{U}\g$ be the Casimir operator with
 respect to the unique invariant scalar product $B:\,\g\times\g\longrightarrow
 \R$ on the Lie algebra $\g$ of the group $G$, which makes the orbit map
 $G\longrightarrow M,\,g\longmapsto g\centerdot p,$ through a point
 $p\,\in\,M$ a Riemannian submersion. The image of $\Cas\,\in\,\mathscr{U}\g$
 under the left regular representation $L$ on $\Gamma(\,\S^2_\circ T^*M\,)$
 satisfies
 \begin{equation}\label{tinv}
  \Delta_L\,\circ\,F
  \;\;=\;\;
  L_\Cas\,\circ\,F
  \;\;=\;\;
  F\,\circ\,(\,\Cas\,\star_R\,)
 \end{equation}
 for every $G$--equivariant linear map $F:\,R\longrightarrow
 \Gamma(\,\S^2_\circ T^*M\,)$. According to equation (\ref{tinv}) every finite
 direct sum $R\,\subseteq\,\Gamma(\,\S^2_\circ T^*M\,)$ of eigenspaces of
 $\Delta_L$ is $G$--characteristic, because every $G$--equivariant linear
 map $F:\,R\longrightarrow\Gamma(\,\S^2_\circ T^*M\,)$ maps $R$ necessarily
 into the direct sum of eigenspaces of $\Delta_L$ for the eigenvalues the
 Casimir operator $\Cas$ assumes $R$:

 \begin{Corollary}[Stability of Einstein Metrics on Symmetric Spaces]
 \hfill\label{lsrs}\break
  Consider a compact irreducible Riemannian symmetric space $M$ of dimension
  $m$ and scalar curvature $\scal$ as a homogeneous space under the transitive
  action $G\times M\longrightarrow M,\,g\longmapsto g\centerdot p$ of its
  isometry group $G$. The direct sum of all eigenspaces of the Lichnerowicz
  Laplacian $\Delta_L$ for the eigenvalues below the critical value
  $2\,\frac\scal m$ is a finite dimensional $G$--characteristic subspace
  $R\,\subseteq\,\Gamma(\,\S^2_\circ T^*M\,)$. In particular the symmetric
  Einstein metric on $M$ is stable
  $$
   R\;\cap\;\ker\;\Big(\;D^*:\;\;\Gamma(\,\S^2_\circ T^*M\,)
   \;\longrightarrow\;\Gamma(\,TM\,)\;\Big)
   \;\;=\;\;
   \{\,0\,\}\ ,
  $$
  if and only if the prototypical differential operator associated to the
  divergence $D^*$ is injective:
  $$
   D^*_R:\;\;\Hom_H(\;R,\;\S^2_\circ T^*_pM\;)\;\longrightarrow\;
   \Hom_H(\;R,\;T_pM\;)\ .
  $$
 \end{Corollary}

 \noindent
 For Riemannian symmetric spaces the critical value $2\,\frac\scal m$ has a
 rather direct interpretation in terms of the Casimir operator $\Cas$. In order
 to establish this interpretation we remark that the standard Laplace operator
 $\Delta$ is defined in \cite{SW19} for every natural vector bundle over a
 Riemannian manifold $M$, not only for $\S^2_\circ T^*M$, in particular the
 version of $\Delta$ defined on the tangent bundle appears in the following
 identity for Killing vector fields $X\,\in\,\Gamma(\,TM\,)$:
 \begin{equation}\label{kvf}
  \Delta X
  \;\;=\;\;
  2\;\Ric\,X\ .
 \end{equation}
 For a Riemannian symmetric space $M$ the fundamental vector field map
 $$
  F:\;\;\g\;\longrightarrow\;\Gamma(\,TM\,),
  \;\;X\;\longmapsto\;X^M\;\;:=\;\;\left.\frac d{dt}\right|_0\mu_{e^{tX}}
 $$
 maps $\g$ to the Lie algebra of Killing vector fields and so equations
 (\ref{tinv}) and (\ref{kvf}) imply:
 $$
  2\;\Ric\,\circ\,F
  \;\;=\;\;
  \Delta\,\circ\,F
  \;\;=\;\;
  F\,\circ\,(\,\Cas\,\star_\g\,)\ .
 $$

 \begin{Remark}[Casimir Eigenvalue of Adjoint Representation]
 \hfill\label{case}\break
  The eigenvalue of the Casimir operator $\Cas\,\in\,\mathscr{U}\g$ on the
  adjoint representation $\g$ of the isometry group $G$ of an irreducible
  Riemannian symmetric space $M$ of dimension $m$ equals
  $$
   \Cas\,\star_\g
   \;\;=\;\;
   2\,\frac\scal m\,\id_\g
  $$
  provided the Casimir operator $\Cas\,\in\,\mathscr{U}\g$ is defined with
  respect to the unique invariant scalar product on $\g$, which makes the
  orbit map $G\longrightarrow M,\,g\longmapsto g\centerdot p,$ a Riemannian
  submersion.
 \end{Remark}

 \noindent
 Corollary \ref{lsrs} and Remark \ref{case} reduce Theorem \ref{mt} in
 essence to a simple and straightforward problem in linear algebra. In order
 to make this reduction effective however we need to calculate the prototypical
 differential operators $D^*_R$ associated to the divergence operator
 $D^*:\,\Gamma(\,\S^2_\circ T^*M\,)\longrightarrow\Gamma(\,TM\,)$ of
 equation (\ref{ll*}) restricted to trace free symmetric $2$--tensors.
 For this purpose we recall that the Lie derivative of a section
 $v\,\in\,\Gamma(\,VM\,)$ of a first order natural vector bundle $VM$
 along $X\,\in\,\Gamma(\,TM\,)$ can be written in terms of the connection
 $$
  \mathfrak{Lie}_Xv
  \;\;=\;\;
  \nabla_Xv\;-\;(\,\nabla^{\mathrm{op}}X\,)\,\star\,v\ ,
 $$
 induced on $VM$ from a connection $\nabla$ on the tangent bundle, where
 $\nabla^{\mathrm{op}}_YX\,:=\,\nabla_XY\,-\,[X,Y]$ denotes the connection
 opposite to $\nabla$ and $\star$ the infinitesimal representation of the
 Lie algebra bundle $\End\,TM$ on $VM$. For the torsion free Levi--Civita
 connection and the fundamental vector field $X^M\,:=\,\left.\frac d{dt}
 \right|_0\mu_{e^{tX}}\,\in\,\Gamma(\,TM\,)$ associated to $X\,\in\,\g$
 this formula becomes:
 \begin{equation}\label{nx}
  \nabla_{X^M}v
  \;\;=\;\;
  \mathfrak{Lie}_{X^M}v\;+\;(\,\nabla X^M\,)\,\star\,v
  \;\;=\;\;
  \left.\frac d{dt}\right|_0\mu_{e^{tX}}^*v
  \;+\;(\,\nabla X^M\,)\,\star\,v\ .
 \end{equation}
 In contrast to general Riemannian homogeneous spaces the Lie algebra
 $\g$ of the isometry group $G$ of a Riemannian symmetric space $M$
 splits naturally into the direct sum $\g\,=\,\h\,\oplus\,\m$ of the
 Lie subalgebra $\h\,\subseteq\,\g$ of the stabilizer $H\,\subseteq\,G$
 of a chosen base point $p\,\in\,M$ and the subspace $\m$ of transvections,
 which correspond to the Killing vector fields which are parallel
 $$
  \m
  \;\;:=\;\;
  \{\;\;X\,\in\,\g\;\;|\;\;(\,\nabla X^M\,)_p\,=\,0\;\;\}
  \;\stackrel\cong\longrightarrow\;T_pM,\qquad X\;\longmapsto\;X^M_p
 $$
 in $p$. Chosing a basis $E_1,\,\ldots,\,E_m\,\in\,\m$ with orthonormal values
 in $T_pM$ we thus find
 \begin{eqnarray*}
  (\,D_R^*F\,)(\,r\,)
  \;\;:=\;\;
  \Big(\;D^*F^\ext(\,r\,)\;\Big)_p
  &=&
  -\;\sum_{\mu\,=\,1}^m\Big(\;E^M_\mu\;\ins\;
  \nabla_{E^M_\mu}F^\ext(\,r\,)\;\Big)^\#_p
  \\
  &=&
  -\;\sum_{\mu\,=\,1}^m\Big(\;(E^M_\mu)^{\;}_p\;\ins\;
  \Big(\;\left.\frac d{dt}\right|_0\mu^*_{e^{tE_\mu}}F^\ext(\,r\,)
  \;\Big)_p\;\Big)^\#
  \\
  &=&
  -\;\sum_{\mu\,=\,1}^m\Big(\;(E^M_\mu)^{\;}_p\;\ins\;
  \Big(\;\left.\frac d{dt}\right|_0L_{e^{-tE_\mu}}F^\ext(\,r\,)
  \;\Big)_p\;\Big)^\#
 \end{eqnarray*}
 for every given argument $F\,\in\,\Hom_H(\,R,\,\S^2_\circ T^*_pM\,)$ with
 Frobenius reciprocal extension $F^\ext:\,R\longrightarrow\Gamma(\,\S^2_\circ
 T^*M\,)$ by using $(\,\nabla E^M_\mu\,)_p\,=\,0$ and equation (\ref{req})
 in the last line. Note that the equality $(\,D^*h\,)_p\,=\,-\sum(\,E^M_\mu\,
 \ins\,\nabla_{E^M_\mu}h\,)^\#_p$ requires orthonormality of the local basis
 $E^M_1,\,\ldots,\,E^M_m$ in the point $p\,\in\,M$ only, because $X\,\ins\,
 \nabla_Yh$ depends $C^\infty(\,M\,)$--bilinearly on $X$ and $Y$. By
 construction $F^\ext$ is $G$--equivariant with $L_{e^{-tE_\mu}}F^\ext
 (\,r\,)\,=\,F^\ext(\,e^{-tE_\mu}\star_Rr\,)$ so that we can summarize our
 arguments for all $F\,\in\,\Hom_H(\,R,\,\S^2_\circ T^*_pM\,)$ in the form:
 \begin{equation}\label{proto}
  (\,D_R^*F\,)(\,r\,)
  \;\;=\;\;
  -\;\sum_{\mu\,=\,1}^m
  \Big(\;(E^M_\mu)^{\;}_p\;\ins\;F(\,-\,E_\mu\,\star_R\,r\,)\;\Big)^\#\ .
 \end{equation}
\section{Critical Representations of Simple Lie Algebras}
\label{critreps}
 On a compact irreducible Riemannian symmetric space $M$ of dimension $m$
 and scalar curvature $\scal$ the direct sum of all eigenspaces of the
 Lichnerowicz Laplacian $\Delta_L$ on trace free symmetric $2$--tensors with
 eigenvalues below the critical value $2\,\frac\scal m$ is a characteristic
 subspace $R\,\subseteq\,\Gamma(\,\S^2_\circ T^*M\,)$ under the isometry
 group $G$. Being a characteristic subspace $R$ equals the direct sum of all
 isotypical components of $\Gamma(\,\S^2_\circ T^*M\,)$ corresponding to
 irreducible representations of $G$ with Casimir eigenvalues below the
 Casimir eigenvalue $\Cas_\g\,=\,2\,\frac\scal m$ of the adjoint
 representation $\g$ according to Remark \ref{case}:
 
 \begin{Definition}[Critical Representation]
 \hfill\label{cr}\break
  An irreducible finite dimensional representation $R$ of a simple Lie
  algebra $\g$ is a critical representation with respect to the stability
  problem of Einstein metrics provided its Casimir eigenvalue $\Cas_R$ is
  at most equal to the Casimir eigenvalue $\Cas_\g$ of the adjoint
  representation:
  $$
   0\;\;\leq\;\;\frac{\Cas_R}{\Cas_\g}\;\;\leq\;\;1\ .
  $$
 \end{Definition}

 \noindent
 In order to classify the critical representations of a simple complex Lie
 algebra $\g$ of rank $n$ we choose a maximal torus $\t\,\subseteq\,\g$, an
 ordering of roots in form of a system $\alpha_1,\,\ldots,\,\alpha_n\,\in\,
 \t^*$ of simple roots and a Weyl invariant scalar product $b:\,\t^*\times\t^*
 \longrightarrow\C$, which is positive definite on the real subspace spanned
 by the (simple) roots. The geometry of the system of simple roots is encoded
 in the Dynkin diagram of the simple Lie algebra $\g$ \cite{fh}, and this
 geometry determines the fundamental weights $\omega_1,\,\ldots,\,
 \omega_n\,\in\,\t^*$ by means of the identity
 $$
  b(\;\omega_r,\;\alpha_s\;)
  \;\;=\;\;
  \frac{\delta_{rs}}2\;b(\;\alpha_s,\;\alpha_s\;)
 $$
 valid for all $r,\,s\,=\,1,\ldots,n$, where $\delta_{rs}$ denotes the
 Kronecker delta which casts the boolean expression $r\,=\,s$ to $1$ or
 $0$. A detailed tabulation of the systems of fundamental weights for all
 simple complex Lie algebras can be found in a number of references, for
 example in \cite{fh}, and we will freely use this information in the
 sequel. Up to isomorphism every irreducible finite dimensional representation
 $R\,=\,R^\lambda$ is characterized by its highest weight, an element
 \begin{equation}\label{hwg}
  \lambda
  \;\;=\;\;
  \lambda_1\,\omega_1\;+\;\ldots\;+\;\lambda_n\,\omega_n
 \end{equation}
 of the additive semigroup of $\t^*$ generated by the fundamental weights
 with $\lambda_1,\,\ldots,\,\lambda_n\,\in\,\N_0$. According to Freudenthal's
 formula for the Casimir eigenvalues \cite{fh} the Casimir operator
 $\Cas\,\in\,\mathscr{U}\g$ can be normalized in such a way that the
 Casimir eigenvalue of every irreducible representation $R\,=\,R^\lambda$
 equals the Casimir constant $\Cas_R\,=\,\Cas_\lambda$ of its highest weight
 $\lambda$:
 \begin{equation}\label{cev}
  \Cas_\lambda
  \;\;:=\;\;
  b(\;\lambda,\;\lambda\,+\,2\,\rho\;)
  \qquad\qquad
  \rho
  \;\;:=\;\;
  \omega_1\;+\;\ldots\;+\;\omega_n\ .
 \end{equation}
 By construction the highest weights of the form (\ref{hwg}) form an additive
 semigroup $\L^+\,\subseteq\,\t^*$. Addition in this semigroup corresponds to
 the Cartan product of irreducible representations, in other words the Cartan
 product of two irreducible representations $R^\lambda$ and $R^{\hat\lambda}$
 of highest weights $\lambda$ and $\hat\lambda$ is some irreducible
 representation $R^{\lambda\,+\,\hat\lambda}\,=\,R^\lambda\odot
 R^{\hat\lambda}$ of highest weight $\lambda\,+\,\hat\lambda$. The Casimir
 eigenvalue of the Cartan product equals the sum of Casimir eigenvalues
 \begin{eqnarray*}
  \Cas_{\lambda+\hat\lambda}
  &=&
  b(\,\lambda+\hat\lambda,\,\lambda+\hat\lambda+2\rho\,)
  \\[2pt]
  &=&
  b(\,\lambda,\,\lambda+2\rho\,)
  \;+\;2\,b(\,\lambda,\,\hat\lambda\,)
  \;+\;b(\,\hat\lambda,\,\hat\lambda+2\rho\,)
  \;\;=\;\;
  \Cas_\lambda\;+\;2\,b(\,\lambda,\,\hat\lambda\,)\;+\;\Cas_{\hat\lambda}
 \end{eqnarray*}
 plus an additional term $b(\,\lambda,\,\hat\lambda\,)$. Both highest weights
 $\lambda$ and $\hat\lambda$ pertain to the Weyl chamber, the real cone in
 $\t^*$ generated by the fundamental weights $\omega_1,\,\ldots,\,\omega_n$,
 and so the additional term is strictly positive $b(\,\lambda,\,\hat\lambda\,)
 \,>\,0$ unless $\lambda\,=\,0$ or $\hat\lambda\,=\,0$. Thought of as a map
 $\Cas:\,\Lambda^+\longrightarrow\R,\,\lambda\longmapsto\Cas_\lambda,$ the
 Casimir eigenvalue is thus superadditive in the sense $\Cas_{\lambda+
 \hat\lambda}\,\geq\,\Cas_\lambda\,+\,\Cas_{\hat\lambda}$, hence the equality
 (\ref{hwg}) converts into the linear lower bound
 \begin{equation}\label{lest}
  \Cas_\lambda
  \;\;\geq\;\;
  \lambda_1\,\Cas_{\omega_1}\;+\;\ldots\;+\;\lambda_n\,\Cas_{\omega_n}
 \end{equation}
 for the Casimir eigenvalue of a highest weight $\lambda\,\in\,\L^+$ with
 equality, if and only if $\lambda$ is itself a fundamental weight or zero.
 This linear lower bound is of course rather weak, after all the Casimir
 eigenvalue is the quadratic polynomial (\ref{cev}) in $\lambda$, nevertheless
 it suffices to select a short list of candidates for critical representations
 for every simple Lie algebra $\g$.

 \pfill
 Let us discuss this procedure for the sequence of special linear Lie
 algebras $\sl(n+1,\C)$ for all $n\,\geq\,1$, which correspond to the
 Dynkin diagrams of type $A_n$. The weights $\e_1,\,\ldots,\,\e_{n+1}$ of
 the defining representation $V\,:=\,\C^{n+1}$ generate the dual $\t^*$
 of the maximal torus $\t\,\subseteq\,\sl(n+1,\C)$ of diagonal matrices
 and are subject only to the characteristic trace zero constraint
 $$
  \e_1\;+\;\e_2\;+\;\ldots\;+\;\e_{n+1}
  \;\;=\;\;
  0\ .
 $$
 which is reflected by the following Weyl invariant scalar product
 $b:\,\t^*\times\t^*\longrightarrow\C$ on $\t^*$:
 \begin{equation}\label{slnpro}
  b(\;\e_\mu,\;\e_\nu\;)
  \;\;:=\;\;
  \delta_{\mu\nu}\;-\;\frac1{n+1}\ .
 \end{equation}
 Due to the trace zero constraint the generating weights $\e_1,\,\ldots,
 \,\e_{n+1}$ are linearly dependent, nevertheless the scalar product
 $b(\,\lambda,\,\hat\lambda\,)$ of two linear combinations $\lambda,
 \,\hat\lambda\,\in\,\t^*$ of $\e_1,\,\ldots,\,\e_{n+1}$ is readily
 calculated by first pretending that $\e_1,\,\ldots,\,\e_{n+1}$ is an
 orthonormal basis of $\t^*$ and then subtracting the product of the
 two coefficient sums of $\lambda$ and $\hat\lambda$ respectively
 divided by $n+1$. For a suitable ordering of roots the simple roots and
 fundamental weights read
 $$
  \begin{array}{lclclcl}
   \alpha_1&=&\e_1\,-\,\e_2&\qquad&
   \omega_1&=&\e_1
   \\
   \alpha_2&=&\e_2\,-\,\e_3&&
   \omega_2&=&\e_1\,+\,\e_2
   \\
   \alpha_3&=&\e_3\,-\,\e_4&&
   \omega_3&=&\e_1\,+\,\e_2\,+\,\e_3
   \\
   &\vdots&&&&\vdots&
   \\
   \alpha_n&=&\e_n\,-\,\e_{n+1}&&
   \omega_n&=&\e_1\,+\,\e_2\,+\,\e_3\,+\,\cdots\,+\,\e_n\ ,
  \end{array}
 $$
 in particular $2\,\rho\,=\,2n\,\e_1\,+\,2(n-1)\,\e_2\,+\,\cdots\,+\,2\,\e_n$
 has coefficient sum equal to $n(n+1)$. The adjoint representation of
 $\g\,=\,\sl(n+1,\C)$ is isomorphic to the Cartan product $V^*\odot V$
 of the irreducible representations for $\omega_n$ and $\omega_1$ of
 highest weight $\e_1-\e_{n+1}$ resulting in:
 $$
  \Cas_\g
  \;\;=\;\;
  b(\,\e_1-\e_{n+1},\,\e_1-\e_{n+1}+2\,\rho\,)
  \;\;=\;\;
  2\,n\,+\,2\,-\,\frac{0\,\cdot\,(\,0+n(n+1)\,)}{n+1}
  \;\;=\;\;
  2\,(n+1)\ .
 $$
 The irreducible representation corresponding to the fundamental weight
 $\omega_r$ equals the exterior power $\L^rV\,\cong\,\L^{n+1-r}V^*$ of the
 defining representation $V$ with Casimir eigenvalue:
 \begin{eqnarray*}
  \Cas_{\omega_r}
  &=&
  b(\,\e_1\,+\,\cdots\,+\,\e_r,\,\e_1\,+\,\cdots\,+\,\e_r\,+\,2\,\rho\,)
  \\
  &=&
  r\,(\,2n-r+2\,)\;-\;\frac{r\,\cdot\,(\,r\,+\,n(n+1)\,)}{n+1}
  \;\;=\;\;
  \frac{n+2}{n+1}\,r\,(\,n+1-r\,)\ .
 \end{eqnarray*}
 Considered as a quadratic polynomial in $r$ the Casimir eigenvalue has zeroes
 in $r\,=\,0,\,n+1$, hence it is symmetric about $r\,=\,\frac{n+1}2$ and
 strictly monotonely increasing for $r\,\in\,[\,0,\,\frac{n+1}2\,]$:
 \begin{equation}\label{simineq}
  \Cas_{\omega_1}
  \;\;<\;\;
  \Cas_{\omega_2}
  \;\;<\;\;
  \Cas_{\omega_3}
  \;\;<\;\;
  \ldots
  \;\;>\;\;
  \Cas_{\omega_{n-2}}
  \;\;>\;\;
  \Cas_{\omega_{n-1}}
  \;\;>\;\;
  \Cas_{\omega_n}\ .
 \end{equation}
 In passing we observe that there are fewer different fundamental weights than
 shown in this diagram for $n\,\leq\,4$. The classification of all possible
 critical representations of the Lie algebras $\sl(n+1,\C)$ with $n\,\geq\,1$
 follows from the lower linear estimate (\ref{lest}), the inequality
 \begin{equation}\label{2most}
  3\,\frac{\Cas_{\omega_1}}{\Cas_\g}
  \;\;=\;\;
  \frac{3\,(n+2)\,n}{2\,(n+1)^2}
  \;\;>\;\;
  1
 \end{equation}
 valid for all $n\,\geq\,1$ and two additional inequalities valid
 for $n\,\geq\,3$ and $n\,\geq\,7$ respectively:
 \begin{equation}\label{56case}
  \frac{\Cas_{\omega_1}}{\Cas_\g}\;+\;\frac{\Cas_{\omega_2}}{\Cas_\g}
  \;\;=\;\;
  \frac{(n+2)\,(3n-2)}{2\,(n+1)^2}
  \;\;>\;\;
  1
  \qquad
  \frac{\Cas_{\omega_3}}{\Cas_\g}
  \;\;=\;\;
  \frac{3\,(n+2)\,(n-2)}{2\,(n+1)^2}
  \;\;>\;\;
  1\ .
 \end{equation}
 According to the inequalities (\ref{simineq}) and (\ref{2most}) the highest
 weight $\lambda$ of every critical representation for $\sl(n+1,\C)$ is the
 sum of at most two fundamental weights, the first inequality in (\ref{56case})
 ensures in turn that $\lambda$ is either zero, a fundamental weight or the
 sum $2\,\omega_1$, $\omega_1+\omega_n$ or $2\,\omega_n$ of two fundamental
 weights in $\{\,\omega_1,\,\omega_n\,\}$. Moreover the second inequality in
 (\ref{56case}) excludes the fundamental weights $\omega_3,\,\ldots,\,
 \omega_{n-2}$ for all $n\,\geq\,7$ from consideration, while the two sums
 $2\,\omega_1$ and $2\,\omega_n$ have to be discarded by hand for all
 $n\,\geq\,2$ by verifying for those $n$:
 $$
  \Cas_{2\,\omega_1}
  \;\;=\;\;
  b(\,2\,\e_1,\,2\,\e_1\,+\,2\,\rho\,)
  \;\;=\;\;
  4\,n\,+\,4\,-\,\frac{2\,\cdot\,(\,2+n(n+1)\,)}{n+1}
  \;\;=\;\;
  \frac{2\,n\,(n+3)}{(n+1)}
  \;\;>\;\;
  \Cas_\g\ .
 $$
 In consequence the trivial representation $\C$, the adjoint representation
 $V^*\odot V$ and the four fundamental representations $V,\,\L^2V,\,\L^2V^*$
 and $V^*$ of $\sl(n+1,\C)$ are critical representations, and they comprise
 all possible critical representations of $\sl(n+1,\C)$ unless $n\,=\,5$ or
 $n\,=\,6$, in which case $\L^3V$ and $\L^3V^*$ are critical representations
 of $\sl(n+1,\C)$ as well.

 \pfill
 Repeating the preceeding discussion for the other three series of classical
 simple Lie algebras we recall first of all that the Dynkin diagrams of type
 $B_n$ with $n\,\geq\,2$ are represented by the odd dimensional orthogonal
 Lie algebras $\so(2n+1,\C)$ of rank $n\,\in\,\N$. The weights $\pm\e_1,\,
 \ldots,\,\pm\e_n$ and $0$ of the defining representation $V\,:=\,\C^{2n+1}$
 of $\so(2n+1,\C)$ form an orthonormal basis $\e_1,\,\ldots,\,\e_n$ of the
 dual $\t^*$ of a maximal torus with respect to the Weyl invariant scalar
 product $b:\,\t^*\times\t^*\longrightarrow\C$ of choice. For a suitable
 ordering of roots the simple roots and fundamental weights of the Lie
 algebras of type $B_n$ read:
 $$
  \begin{array}{lclclcl}
   \alpha_1&=&\e_1\,-\,\e_2&\qquad&
   \omega_1&=&\phantom{\frac12(\,}\e_1
   \\
   \alpha_2&=&\e_2\,-\,\e_3&&
   \omega_2&=&\phantom{\frac12(\,}\e_1\,+\,\e_2
   \\
   &\vdots&&&&\vdots&
   \\
   \alpha_{n-1}&=&\e_{n-1}\,-\,\e_n&&
   \omega_{n-1}&=&\phantom{\frac12(\,}\e_1\,+\,\e_2\,+\,\cdots\,+\,\e_{n-1}
   \\
   \alpha_n&=&\e_n&&
   \omega_n&=&\frac12(\,\e_1\,+\,\e_2\,+\,\cdots\,+\,\e_{n-1}\,+\,\e_n\,)\ .
  \end{array}
 $$
 Using $2\,\rho\,=\,(2n-1)\,\e_1\,+\,(2n-3)\,\e_2\,+\,\cdots\,+\,\e_n$ we
 can calculate the Casimir eigenvalues
 $$
  \Cas_{\omega_r}
  \;\;=\;\;
  b(\;\e_1\,+\,\cdots\,+\,\e_r,\;\e_1\,+\,\cdots\,+\,\e_r\,+\,2\,\rho\;)
  \;\;=\;\;
  r\,(\,2n\,+\,1\,-\,r\,)
 $$
 for the exterior powers $\L^rV$ of the defining representation $V$, which
 are irreducible representations of highest weight $\omega_r$ for $r\,=\,
 1,\,\ldots,\,n-1$. In the same vein the Casimir eigenvalue of the spinor
 representation, the irreducible representation $\Z$ of highest weight
 $\omega_n$, equals:
 $$
  \Cas_{\omega_n}
  \;\;=\;\;
  b(\;{\textstyle\frac12}\,(\,e_1\,+\,\cdots\,+\,\e_n\,),\;
  {\textstyle\frac12}\,(\,\e_1\,+\,\cdots\,+\,\e_n\,)\,+\,2\,\rho\;)
  \;\;=\;\;
  {\textstyle\frac14}\,n\,(\,2n+1\,)\ .
 $$
 The Casimir eigenvalue of the adjoint representation $\L^2V$ of $\so(2n+1,\C)$
 equals $2\,(2n-1)$ even in the exceptional case $n\,=\,2$, where the highest
 weight of $\L^2V$ equals $2\,\omega_2\,=\,\e_1+\e_2$. Discussing the roots of
 the quadratic polynomial $\Cas_{\omega_r}\,=\,r\,(2n+1-r)$ in $r$ we conclude:
 \begin{equation}\label{bnseq}
  \Cas_{\omega_1}
  \;\;<\;\;
  \Cas_{\omega_2}
  \;\;<\;\;
  \ldots
  \;\;<\;\;
  \Cas_{\omega_{n-1}}\ .
 \end{equation}
 Taking these inequalities together with the lower bound (\ref{lest})
 and the auxiliary inequalities
 $$
  \frac{\Cas_{\omega_1}}{\Cas_\g}
  \;\;=\;\;
  \frac{2\,n}{2\,(2n-1)}
  \;\;>\;\;
  \frac12
  \qquad\qquad
  \frac{\Cas_{\omega_n}}{\Cas_\g}
  \;\;=\;\;
  \frac{2n\,(2n+1)}{16\,(2n-1)}
  \;\;>\;\;
  \frac12
 $$
 valid for all $n\,\geq\,2$ and $n\,\geq\,3$ respectively we conclude that
 the highest weight $\lambda\,\in\,\t^*$ of a critical representation for
 $\so(2n+1,\C)$ is either zero or a fundamental weight for all $n\,\geq\,3$,
 interestingly $2\,\omega_2$ is the highest weight of the critical adjoint
 representation $\so(5,\C)\,=\,\L^2V$ in the exceptional case $n\,=\,2$.
 Due to the inequalities (\ref{bnseq}) combined with $\Cas_{\omega_2}\,=\,
 \Cas_\g$ the fundamental weights $\omega_3,\,\ldots,\,\omega_{n-1}$ can not
 be critical for $n\,\geq\,3$, in consequence classification of the critical
 irreducible representations of $\so(2n+1,\C)$ centers on the candidates
 $\C,\,V,\,\L^2V$ and the spinor representation $\Sigma$, the latter however
 is critical only for $n\,=\,2,\,\ldots,\,6$.

 \pfill
 Let us now turn to the symplectic Lie algebras $\g\,=\,\sp(2n,\C)$ of
 rank $n\,\in\,\N$ and Dynkin diagrams of type $C_n$ under the assumption
 $n\,\geq\,2$. The weights $\pm\e_1,\,\ldots,\,\pm\e_n$ of the defining
 representation $V\,=\,\C^{2n}$ of the symplectic Lie algebras form a basis
 $\e_1,\,\ldots,\,\e_n$ of the dual $\t^*$ of a maximal torus $\t$, which
 is actually an orthonormal basis for the Weyl--invariant scalar product
 $b:\,\t^*\times\t^*\longrightarrow\C$ of our choice. In terms of this
 orthonormal basis the simple roots associated to a suitable ordering of
 roots and the corresponding fundamental weights equal
 $$
  \begin{array}{lclclcl}
   \alpha_1&=&\e_1\,-\,\e_2&\qquad&
   \omega_1&=&\e_1
   \\
   \alpha_2&=&\e_2\,-\,\e_3&&
   \omega_2&=&\e_1\,+\,\e_2
   \\
   &\vdots&&&&\vdots&
   \\
   \alpha_{n-1}&=&\e_{n-1}\,-\,\e_n&&
   \omega_{n-1}&=&\e_1\,+\,\e_2\,+\,\cdots\,+\,\e_{n-1}
   \\
   \alpha_n&=&2\,\e_n&&
   \omega_n&=&\e_1\,+\,\e_2\,+\,\cdots\,+\,\e_{n-1}\,+\,\e_n\ ,
  \end{array}
 $$
 and so $2\,\rho\,=\,2n\,\e_1\,+\,2(n-1)\,\e_2\,+\,\cdots\,+\,2\,\e_n$.
 In consequence the Casimir eigenvalue of the fundamental representation
 $\L^r_\circ V$ corresponding to the fundamental weight $\omega_r$ equals
 $$
  \Cas_{\omega_r}
  \;\;=\;\;
  b(\;\e_1\,+\,\cdots\,+\,\e_r,\;\e_1\,+\,\cdots\,+\,\e_r\,+\,2\,\rho\;)
  \;\;=\;\;
  r\,(\,2n\,-\,r\,+\,2\,)
 $$
 for $r\,=\,1,\,\ldots,\,n$. In analogy to the inequalities (\ref{simineq})
 and (\ref{bnseq}) we establish the inequalities:
 \begin{equation}\label{cnseq}
  \Cas_{\omega_1}
  \;\;<\;\;
  \Cas_{\omega_1}
  \;\;<\;\;
  \ldots
  \;\;<\;\;
  \Cas_{\omega_n}\ .
 \end{equation}
 The adjoint representation of $\g\,=\,\sp(2n,\C)$ is isomorphic to the
 second symmetric power $\sp(2n,\C)\,=\,\S^2V$ of the defining representation
 $V$, its Casimir eigenvalue equals:
 $$
  \Cas_{2\,\e_1}
  \;\;=\;\;
  b(\;2\,\e_1,\;2\,\e_1\,+\,2\,\rho\;)
  \;\;=\;\;
  4\,(n+1)\ .
 $$
 The sequence of inequalities (\ref{cnseq}) together with the auxiliary
 inequalities
 $$
  \frac{\Cas_{\omega_1}}{\Cas_\g}\;+\;\frac{\Cas_{\omega_2}}{\Cas_\g}
  \;\;=\;\;
  \frac{6n+1}{4\,(n+1)}
  \;\;>\;\;
  1
  \qquad\qquad
  \frac{\Cas_{\omega_3}}{\Cas_\g}
  \;\;=\;\;
  \frac{3\,(2n-1)}{4\,(n+1)}
  \;\;>\;\;
  1
 $$
 valid for all $n\,\geq\,2$ and all $n\,\geq\,4$ respectively tell us
 that the highest weight $\lambda$ of a critical representation of
 $\sp(2n,\C)$ is either zero $\lambda\,=\,0$, the highest weight
 $\lambda\,=\,2\,\omega_1$ of the adjoint representation or one of the
 two fundamental weights $\lambda\,=\,\omega_1$ or $\lambda\,=\,\omega_2$
 unless $n\,=\,3$, where $\lambda\,=\,\omega_3$ is an additional possibility.
 The critical representations of the symplectic Lie algebras $\sp(2n,\C)$ are
 thus $\C,\,V,\,\S^2V$ and $\L^2_\circ V$ with $\L^3_\circ V$ being critical
 only for $n\,=\,3$.

 \pfill
 The last sequence $D_n$ of irreducible Dynkin diagrams is represented
 by the even dimensional orthogonal Lie algebras $\so(2n,\C)$ of rank
 $n\,\geq\,3$, the classification of their critical representation thus
 follows the discussion of the odd dimensional orthogonal Lie algebras
 closely. The weights $\pm\e_1,\,\ldots,\,\pm\e_n$ of the defining
 representation $V\,:=\,\C^{2n}$ form again an orthonormal basis
 $\e_1,\,\ldots,\,\e_n$ of the dual $\t^*$ of a maximal torus with
 respect to our preferred Weyl invariant scalar product $b:\,\t^*\times
 \t^*\longrightarrow\C$. Simple roots and fundamental weights read
 $$
  \begin{array}{lclclcl}
   \alpha_1&=&\e_1\,-\,\e_2&\qquad&
   \omega_1&=&\phantom{\frac12(\,}\e_1
   \\
   \alpha_2&=&\e_2\,-\,\e_3&&
   \omega_2&=&\phantom{\frac12(\,}\e_1\,+\,\e_2
   \\
   &\vdots&&&&\vdots&
   \\
   \alpha_{n-2}&=&\e_{n-2}\,-\,\e_{n-1}&&
   \omega_{n-2}&=&\phantom{\frac12(\,}\e_1\,+\,\e_2\,+\,\cdots\,+\,\e_{n-2}
   \\
   \alpha_{n-1}&=&\e_{n-1}\,-\,\e_n&&
   \omega_{n-1}&=&\frac12(\,\e_1\,+\,\e_2\,+\,\cdots\,+\,\e_{n-2}
   \,+\,\e_{n-1}\,-\,\e_n\,)
   \\
   \alpha_n&=&\e_{n-1}\,+\,\e_n&&
   \omega_n&=&\frac12(\,\e_1\,+\,\e_2\,+\,\cdots\,+\,\e_{n-2}
   \,+\,\e_{n-1}\,+\,\e_n\,)
  \end{array}
 $$
 for a suitable ordering of roots. The Casimir eigenvalues for the exterior
 powers $\L^rV$ equal
 $$
  \Cas_{\omega_r}
  \;\;=\;\;
  b(\;\e_1\,+\,\cdots\,+\,\e_r,\;\e_1\,+\,\cdots\,+\,\e_r\,+\,2\,\rho\;)
  \;\;=\;\;
  r\,(\,2n\,-\,r\,)\ .
 $$
 for all $r\,=\,1,\,\ldots,\,n-2$ due to $2\,\rho\,=\,(2n-2)\,\e_1\,+\,
 (2n-4)\,\e_2\,+\,\cdots\,+\,2\,\e_{n-1}$. The Casimir eigenvalues of the
 two half spinor representations $\Z^-$ and $\Z^+$ of $\so(2n,\C)$ of highest
 weight $\omega_{n-1}$ and $\omega_n$ agree due to the existence of the
 exterior automorphism $+\e_n\,\rightsquigarrow\,-\e_n$:
 $$
  \Cas_{\omega_{n-1}}
  \;\;=\;\;
  \Cas_{\omega_n}
  \;\;=\;\;
  b(\;{\textstyle\frac12}\,(\,\e_1\,+\,\cdots\,+\,\e_n\,),\;
  {\textstyle\frac12}\,(\,\e_1\,+\,\cdots\,+\,\e_n\,)\,+\,2\,\rho\;)
  \;\;=\;\;
  {\textstyle\frac14}\,n\,(\,2n-1\,)\ .
 $$
 In consequence the Casimir eigenvalue of the adjoint representation $\L^2V$
 of $\so(2n,\C)$ equals $2\,(2n-2)$ even in the exceptional case $n\,=\,3$,
 where the highest weight of $\L^2V$ equals $\omega_2+\omega_3\,=\,\e_1+\e_2$.
 The sequence of inequalities between the Casimir eigenvalues
 \begin{equation}\label{dnseq}
  \Cas_{\omega_1}
  \;\;<\;\;
  \Cas_{\omega_2}
  \;\;<\;\;
  \ldots
  \;\;<\;\;
  \Cas_{\omega_{n-2}}
 \end{equation}
 established by considering $\Cas_{\omega_r}\,=\,r(2n-r)$ as a quadratic
 polynomial in $r$ together with
 $$
  \frac{\Cas_{\omega_1}}{\Cas_\g}
  \;\;=\;\;
  \frac{2n-1}{2\,(2n-2)}
  \;\;>\;\;\frac12
  \qquad\qquad
  \frac{\Cas_{\omega_{n-1}}}{\Cas_\g}
  \;\;=\;\;
  \frac{\Cas_{\omega_n}}{\Cas_\g}
  \;\;=\;\;
  \frac{2n\,(2n-1)}{16\,(2n-2)}
  \;\;>\;\;\frac12
 $$
 for all $n\,\geq\,3$ and $n\,\geq\,4$ respectively imply that the highest
 weight $\lambda$ of a critical representation is either zero or a fundamental
 weight different from $\omega_3,\,\ldots,\,\omega_{n-2}$ unless $n\,=\,3$,
 in which case $\omega_2+\omega_3\,=\,\e_1+\e_2$ equals the highest weight
 of the critical adjoint representation $\L^2V$ of $\so(6,\C)\,\cong\,
 \sl(4,\C)$. Leaving aside this special case already discussed above we
 conclude that the only candidates for critical representations are
 $\C,\,V,\,\L^2V$ and the two half spinor representations $\Z^+$ and
 $\Z^-$, the latter however are critical only for $n\,=\,3,\,\ldots,\,7$.

 \pfill
 For the convenience of the reader we summarize the preceeding classification
 of the critical representations of the classical simple Lie algebras of types
 $A$ to $D$ in the following table, in which the column denoted by Casimir is
 reserved for the relative Casimir eigenvalue $\frac{\Cas_\lambda}{\Cas_\g}$:

 \begin{center}
  \begin{tabular}{| c | c | c | c | l |}
   \hline&&&&
   \\[-10pt]
   \ \small{Algebra}\ & \ \small{Weight}\  & \ \small{Representation}\  &
   \ \small{Casimir}\ & \quad\small{Constraints}
   \\[3pt]
   \hline\hline&&&&
   \\[-10pt]
   $\sl(\,n+1,\,\C\,)$ & $0$ & $\C$ & $0$ &
   \quad for all $n\,\geq\,1$
   \\[3pt]
   \cline{2-5}&&&&
   \\[-10pt]
   & $\omega_1$ & $V$ & $\frac{n(n+2)}{2(n+1)^2}$ &
   \quad for all $n\,\geq\,1$
   \\[3pt]
   \cline{2-5}&&&&
   \\[-10pt]
   & $\omega_n$ & $V^*$ & $\frac{n(n+2)}{2(n+1)^2}$ &
   \quad for all $n\,\geq\,2$
   \\[3pt]
   \cline{2-5}&&&&
   \\[-10pt]
   & $\omega_2$ & $\L^2V$ & $\frac{(n-1)(n+2)}{(n+1)^2}$ &
   \quad for all $n\,\geq\,3$
   \\[3pt]
   \cline{2-5}&&&&
   \\[-10pt]
   & $\omega_{n-1}$ & $\L^2V^*$ & $\frac{(n-1)(n+2)}{(n+1)^2}$ &
   \quad for all $n\,\geq\,4$
   \\[3pt]
   \cline{2-5}&&&&
   \\[-10pt]
   & $\omega_3$ & $\L^3V$ & $\frac78,\;\frac{48}{49}$ &
   \quad for $n\,=\,5,\,6$
   \\[3pt]
   \cline{2-5}&&&&
   \\[-10pt]
   & $\omega_4$ & $\L^3V^*$ & $\frac{48}{49}$ &
   \quad for $n\,=\,6$
   \\[3pt]
   \cline{2-5}&&&&
   \\[-10pt]
   & $\omega_1+\omega_n$ & $V^*\odot V$ & $1$ &
   \quad for all $n\,\geq\,1$
   \\[3pt]
   \hline\hline&&&&
   \\[-10pt]
   $\so(\,2n+1,\,\C\,)$ & $0$ & $\C$ & $0$ &
   \quad for all $n\,\geq\,1$
   \\[3pt]
   \cline{2-5}&&&&
   \\[-10pt]
   & $\omega_1$ & $V$ & $\frac{2n}{2(2n-1)}$ &
   \quad for all $n\,\geq\,1$
   \\[3pt]
   \cline{2-5}&&&&
   \\[-10pt]
   & $\omega_2$ & $\L^2V$ & $1$ &
   \quad for all $n\,\geq\,2$
   \\[3pt]
   \cline{2-5}&&&&
   \\[-10pt]
   & $\omega_n$ & $\Z$ & $\frac{2n(2n+1)}{16(2n-1)}$ &
   \quad for $n\,=\,1,\,\ldots,\,6$\quad\
   \\[3pt]
   \hline\hline&&&&
   \\[-10pt]
   $\sp(\,2n,\,\C\,)$ & $0$ & $\C$ & $0$ &
   \quad for all $n\,\geq\,1$
   \\[3pt]
   \cline{2-5}&&&&
   \\[-10pt]
   & $\omega_1$ & $V$ & $\frac{2n+1}{4(n+1)}$ &
   \quad for all $n\,\geq\,1$
   \\[3pt]
   \cline{2-5}&&&&
   \\[-10pt]
   & $\omega_2$ & $\L^2_\circ V$ & $\frac{4n}{4(n+1)}$ &
   \quad for all $n\,\geq\,2$
   \\[3pt]
   \cline{2-5}&&&&
   \\[-10pt]
   & $\omega_3$ & $\L^3_\circ V$ & $\frac{15}{16}$ &
   \quad for $n\,=\,3$\quad\ 
   \\[3pt]
   \cline{2-5}&&&&
   \\[-10pt]
   & $2\,\omega_1$ & $\S^2V$ & $1$ &
   \quad for all $n\,\geq\,1$
   \\[3pt]
   \hline\hline&&&&
   \\[-10pt]
   $\so(\,2n,\,\C\,)$ & $0$ & $\C$ & $0$ &
   \quad for all $n\,\geq\,3$
   \\[3pt]
   \cline{2-5}&&&&
   \\[-10pt]
   & $\omega_1$ & $V$ & $\frac{2n-1}{2(2n-2)}$ &
   \quad for all $n\,\geq\,3$
   \\[3pt]
   \cline{2-5}&&&&
   \\[-10pt]
   & $\omega_2$ & $\L^2V$ & $1$ &
   \quad for all $n\,\geq\,3$
   \\[3pt]
   \cline{2-5}&&&&
   \\[-10pt]
   & $\omega_{n-1} $ & $\Z^-$ & $\frac{2n(2n-1)}{16(2n-2)}$ &
   \quad for $n\,=\,3,\,\ldots,\,7$\quad\ 
   \\[3pt]
   \cline{2-5}&&&&
   \\[-10pt]
   & $\omega_n $ & $\Z^+$ & $\frac{2n(2n-1)}{16(2n-2)}$ &
   \quad for $n\,=\,3,\,\ldots,\,7$\quad\ 
   \\[3pt]
   \hline
  \end{tabular}
 \end{center}
 The classification of the critical representations of the exceptional simple
 Lie algebras of types $E_6,\,E_7,\,E_8$ and $F_4,\,G_2$ is significantly
 simpler, because all the Casimir eigenvalues are explicit rational numbers
 and the highest weights of the adjoint representations $\g$ are all
 fundamental weights. Moreover all Casimir eigenvalues of the fundamental
 weights turn out to be greater than or equal to $\frac12$ so that the highest
 weights of all non--trivial critical representations are necessarily
 fundamental weights. The resulting classification reads:
 \begin{center}
  \begin{tabular}[t]{| c | c | c |}
   \hline&&
   \\[-10pt]
   \ \small{Algebra}\ & \ \small{Representation}\  & \ \small{Casimir}\ 
   \\[3pt]
   \hline\hline&&
   \\[-10pt]
   $\mathfrak{e}_6$ & $\C$ & $0$
   \\[3pt]
   \cline{2-3}&&
   \\[-10pt]
   & $[\,27\,]$ & $\frac{13}{18}$
   \\[3pt]
   \cline{2-3}&&
   \\[-10pt]
   & $[\,27\,]\hbox to0pt{$^*$\hss}$ & $\frac{13}{18}$
   \\[3pt]
   \cline{2-3}&&
   \\[-10pt]
   & $[\,78\,]$ & $1$
   \\[3pt]
   \hline\hline&&
   \\[-10pt]
   $\mathfrak{e}_7$ & $\C$ & $0$
   \\[3pt]
   \cline{2-3}&&
   \\[-10pt]
   & $[\,56\,]$ & $\frac{19}{24}$
   \\[3pt]
   \cline{2-3}&&
   \\[-10pt]
   & $[\,133\,]$ & $1$
   \\[3pt]
   \hline
  \end{tabular}
  \qquad\qquad
  \begin{tabular}[t]{| c | c | c |}
   \hline&&
   \\[-10pt]
   \ \small{Algebra}\ & \ \small{Representation}\  & \ \small{Casimir}\ 
   \\[3pt]
   \hline\hline&&
   \\[-10pt]
   $\mathfrak{e}_8$ & $\C$ & $0$
   \\[3pt]
   \cline{2-3}&&
   \\[-10pt]
   & $[\,248\,]$ & $1$
   \\[3pt]
   \hline\hline&&
   \\[-10pt]
   $\mathfrak{f}_4$ & $\C$ & $0$
   \\[3pt]
   \cline{2-3}&&
   \\[-10pt]
   & $\mathrm{Im}\;\A$ & $\frac23$
   \\[3pt]
   \cline{2-3}&&
   \\[-10pt]
   & $[\,52\,]$ & $1$
   \\[3pt]
   \hline\hline&&
   \\[-10pt]
   $\mathfrak{g}_2$ & $\C$ & $0$
   \\[3pt]
   \cline{2-3}&&
   \\[-10pt]
   & $\mathrm{Im}\;\O$ & $\frac12$
   \\[3pt]
   \cline{2-3}&&
   \\[-10pt]
   & $[\,14\,]$ & $1$
   \\[3pt]
   \hline
  \end{tabular}
 \end{center}
 In this table we have indicated the non--trivial critical irreducible
 representations by their dimensions only in order to avoid any ambiguity
 caused by the rather arbitrary enumeration of simple roots and fundamental
 weights of the exceptional Lie algebras. The two exceptions are the two
 critical representations of $\mathbf{G}_2$ and $\F_4$ arising from the
 octonions $\O$ and the $27$--dimensional Albert algebra $\A$, which will
 be studied in detail in Section \ref{geocay}.
\section{Geometry of the Cayley Projective Plane}
\label{geocay}
 In order to prove the stability of the Cayley projective plane $\O P^2$
 we will study the Albert algebra $\A$ of hermitean $3\times3$--matrices
 over the octonions $\O$ to describe the unique non--trivial critical
 representation $R\,=\,\Im\,\A$ of its automorphism group $\F_4\,=\,\Aut\,\A$.
 Using the machinery of prototypical differential operators introduced in
 Section \ref{prodiff} and our description of $R$ we will proceed to prove
 the injectivity of the associated prototypical divergence operator $D^*_R$
 and conclude that the Einstein metric on $\O P^2$ is stable in the sense
 of Koiso. A detailed introduction to the octonions and the Albert algebra
 can be found in \cite{baez} or the book \cite{harvey}.

 \pfill
 Working with the algebra $\O$ of Cayley numbers or octonions is hampered
 by the fact that many different constructions and definitions of $\O$
 exist in the literature, it can hardly be called obvious that all these
 constructions lead to the same algebra. Perhaps the most direct construction
 of $\O$ is via the Cayley--Dickson process, in this construction octonions
 are tuples $(\,a,\,\alpha\,)\,\in\,\H\oplus\H$ of quaternions with the
 multiplication, conjugation and algebra unit:
 $$
  (\,a,\,\alpha\,)\,(\,b,\,\beta\,)
  \;\;:=\;\;
  (\,ab\,-\,\overline\beta\alpha,\,\beta a\,+\,\alpha\overline b\,)
  \qquad\quad
  \overline{(\,a,\,\alpha\,)}
  \;\;:=\;\;
  (\,\overline a,-\alpha\,)
  \qquad\quad
  \1
  \;\;:=\;\;
  (\,1,\,0\,)\ .
 $$
 All our subsequent calculations are independent from the question, which
 construction of the octonions is to be preferred; they only rely on four
 characteristic properties of the algebra $\O$, which the reader may easily
 verify in her or his favorite model for the octonions:
 \begin{itemize}
  \item The algebra $\O$ of octonions is an algebra with unit $\1\,\in\,\O$
        which decomposes into the vector space direct sum $\O\,=\,\R\,\oplus
        \,\mathrm{Im}\,\O$ of the line generated by $\1$ and the subspace
        $\mathrm{Im}\,\O\,\subseteq\,\O$ of imaginary octonions spanned by
        all the square roots of $-\1$ in $\O$.
  \item The involutive linear map $\O\longrightarrow\O,\,A\longmapsto
        \overline{A},$ with eigenspaces $\R$ and $\mathrm{Im}\,\O$
        for the eigenvalues $+1$ and $-1$ respectively is an algebra
        antiautomorphism $\overline{A\,B}\,=\,\overline B\,\overline A$.
  \item Every two octonions $A,\,B\,\in\,\O$ lie in a common associative
        subalgebra of $\O$ together with the unit $\1\,\in\,\O$ and in
        consequence $\overline A\,=\,2(\Re\,A)\1\,-\,A$ as well as
        $\overline B$. In particular there is no need to indicate
        parentheses in expressions like $A\overline AB$ or $BA\overline{A}$.
  \item The expression $g_\O(A,B)\,:=\,\Re(\overline A\,B)$ defines a positive
        definite scalar product on $\O$. In the same vein $\Omega(A,B,C)
        \,:=\,\Re(ABC)$ is well--defined for all $A,\,B,\,C\,\in\,\O$
        independent of the way we evaluate it and defines a cyclically
        invariant $3$--form $\Omega$.
 \end{itemize}
 In the Cayley--Dickson construction of the octonions for example the $3$--form
 $\Omega$ reads
 $$
  \Omega\big(\;(\,a,\,\alpha\,),\;(\,b,\,\beta\,),\;(\,c,\,\gamma\,)\;\big)
  \;\;:=\;\;
  \Re\,\big(\;abc\;-\;\overline\alpha\gamma b
  \;-\;\overline\beta\alpha c\;-\;\overline\gamma\beta a\;\big)
 $$
 and is thus cyclically invariant, it becomes an alternating $3$--form
 only after restriction to the subspace $\mathrm{Im}\,\O\,\subseteq\,\O$
 of imaginary octonions. The four properties formulated above in particular
 imply the fundamental identity $\overline{A}A\,=\,|\,A\,|^2_\O\,\1$ for all
 $A\,\in\,\O$, which is needed in the verification of essentially all
 the formulas stipulated below. The Albert algebra is the commutative
 algebra $\A$ of all hermitean $3\times3$--matrices with coefficients
 in $\O$
 $$
  \A
  \;\;:=\;\;
  \left\{\left.\;\;\pmatrix{a_1&\overline A_3&A_2\cr A_3&a_2&\overline A_1\cr
   \overline A_2&A_1&a_3}\;\;\right|\;\;a_1,\,a_2,\,a_3\,\in\,\R
  \textrm{\ and\ }A_1,\,A_2,\,A_3\,\in\,\O\;\;\right\}
 $$
 under the symmetrized matrix multiplication $\mathfrak{A}\,*\,
 \hat\mathfrak{A}\,=\,\frac12(\,\mathfrak{A}\hat\mathfrak{A}\,+\,
 \hat\mathfrak{A}\mathfrak{A}\,)$ or equivalently under:
 $$
  \pmatrix{a_1 & \overline A_3 & A_2\cr A_3 & a_2 & \overline A_1
  \cr\overline A_2 & A_1 & a_3}^2
  =\,
  \pmatrix{
   a_1^2+|A_2|^2+|A_3|^2
   & (a_1+a_2)\overline A_3+A_2A_1
   & (a_3+a_1)A_2+\overline A_3\overline A_1
   \cr
   (a_3+a_1)A_3+\overline A_1\overline A_2
   & a_2^2+|A_3|^2+|A_1|^2
   & (a_2+a_3)\overline A_1+A_3A_2
   \cr
   (a_3+a_1)\overline A_2+A_1A_3
   & (a_2+a_3)A_1+\overline A_2\overline A_3
   & a_3^2+|A_1|^2+|A_2|^2}\ .
 $$
 Let us point out two subspaces of the Albert algebra, the vector and the
 spinor subspace
 \begin{equation}\label{asubs}
  V
  \;\;:=\;\;
  \left\{\;\;\pmatrix{0 & 0 & 0 \cr 0 & +A_0 &
   \hphantom{+}\overline A_1 \cr 0 & \hphantom{+}A_1 & -A_0}\;\;\right\}
  \qquad\qquad
  \Z
  \;\;:=\;\;
  \left\{\;\;\pmatrix{0 & \overline\alpha_3 & \alpha_2 \cr \alpha_3 & 0 & 0
   \cr\overline\alpha_2 & 0 & 0}\;\;\right\}
 \end{equation}
 with parameters $A_0\,\in\,\R$ and $A_1,\,\alpha_2,\,\alpha_3\,\in\,\O$.
 The vector space underlying the Algebra algebra $\A$ decomposes into the
 internal direct sum of $V$ and $\Z$ as well as the two real lines
 \begin{equation}\label{adec}
  \A
  \;\;=\;\;
  \R\;\oplus\;\R\,\Gamma\;\oplus\;V\;\oplus\;\Z\ ,
 \end{equation}
 spanned by the algebra unit $\1\,\in\,\A$ and the auxiliary trace free
 element $\Gamma\,\in\,\A$ defined by:
 \begin{equation}\label{2el}
  \1
  \;\;:=\;\;
  \pmatrix{+1 & 0 & 0 \cr 0 & +1 & 0 \cr 0 & 0 & +1}
  \qquad\qquad
  \Gamma
  \;\;:=\;\;
  \pmatrix{+2 & 0 & 0 \cr 0 & -1 & 0 \cr 0 & 0 & -1}\ .
 \end{equation}
 Most importantly the multiplication in the Albert algebra can be written
 completely in terms of the decomposition (\ref{adec}) and the spin geometry
 of the euclidian vector space $V$, namely
 \begin{eqnarray}
  \Big(\;a\;\oplus\;\underline a\,\Gamma\;\oplus\;A\;\oplus\;\alpha\;\Big)^2
  &=&
  \Big(\;a^2\;+\;2\,\underline a^2\;+\;\frac23\,g_V(\,A,\,A\,)
  \;+\;\frac23\,g_\Z(\,\alpha,\,\alpha\,)\;\Big)
  \nonumber
  \\
  &&
  \;\oplus\;\Big(\;2\,a\,\underline a\;+\;\underline a^2\;-\;
  \frac13\,g_V(\,A,\,A\,)\;+\;\frac16\,g_\Z(\,\alpha,\,\alpha\,)\;\Big)
  \;\Gamma
  \label{amult}
  \\
  &&
  \;\oplus\;\Big(\;(\,2a-2\underline a\,)\,A
  \;+\;\frac12\,\alpha\,\diamond\,\alpha\;\Big)\;\oplus\;
  \Big(\;(\,2a+\underline a\,)\,\alpha\;+\;A\,\bullet\,\alpha\;\Big)
  \nonumber
 \end{eqnarray}
 for all $a,\,\underline a\,\in\,\R$, $A\,\in\,V$ and $\alpha\,\in\,\Z$.
 In this spinorial description of the multiplication in the Albert algebra
 the notation $g_V:\,V\times V\longrightarrow\R$ and $g_\Z:\,\Z\times\Z
 \longrightarrow\R$ refers to the standard positive definite scalar products
 on $V\,\cong\,\R\,\oplus\,\O$ and $\Z\,\cong\,\O\,\oplus\,\O$ respectively
 $$
  g_V(\,A_0\oplus A_1,\,A_0\oplus A_1\,)
  \;\;:=\;\;
  A_0^2\,+\,|\,A_1\,|^2_\O
  \qquad
  g_\Z(\,\alpha_2\oplus\alpha_3,\,\alpha_2\oplus\alpha_3\,)
  \;\;:=\;\;
  |\,\alpha_2\,|^2_\O\,+\,|\,\alpha_3\,|^2_\O\ ,
 $$
 while $\bullet:\,V\times\Z\longrightarrow\Z$ refers to the Clifford
 multiplication of $V$ on its spinor module $\Z$:
 $$
  (\,A_0\,\oplus\,A_1\,)\;\bullet\;(\,\alpha_2\,\oplus\,\alpha_3\,)
  \;\;:=\;\;
  \Big(\;-\,A_0\,\alpha_2\;+\;\overline\alpha_3\,\overline A_1\;\Big)
  \;\oplus\;
  \Big(\;+\,A_0\,\alpha_3\;+\;\overline A_1\,\overline\alpha_2\;\Big)\ .
 $$
 Last but not least $\diamond:\,\Z\times\Z\longrightarrow V$ denotes the
 symmetric bilinear spinor multiplication:
 $$
  (\,\alpha_2\,\oplus\,\alpha_3\,)\;\diamond\;(\,\beta_2\,\oplus\,\beta_3\,)
  \;\;:=\;\;
  \Big(\;g_\O(\alpha_3,\beta_3)\;-\;g_\O(\alpha_2,\beta_2)\;\Big)
  \;\oplus\;
  \Big(\;\overline\alpha_2\,\overline\beta_3\;+\;\overline\beta_2\,
  \overline\alpha_3\;\Big)\ .
 $$
 The verification of the formula (\ref{amult}) for the multiplication in
 $\A$ is a slightly tedious, but otherwise straightforward exercise in
 expanding definitions. In the same vein the identities
 \begin{equation}\label{2id}
  A\bullet(\,A\bullet\alpha\,)
  \;\;=\;\;
  g_V(A,A)\,\alpha
  \qquad\qquad
  (\,\alpha\diamond\alpha\,)\bullet\alpha
  \;\;=\;\;
  g_\Z(\alpha,\alpha)\,\alpha
 \end{equation}
 are easily verified for all $A\,\in\,V$ and $\alpha\,\in\,\Z$ by
 expanding definitions and using the cyclic invariance of $\Omega$.
 It is the first of these two identities of course which allows us to think
 of $\bullet$ as the multiplication of the Clifford algebra $\Cl(\,V,\,
 -g_V\,)\,\cong\,\mathrm{Mat}_{16\times16}(\R\,\oplus\,\R)$ on its
 irreducible module $\Z\,\cong\,\R^{16}$. The cyclic invariance of
 $\Omega$ is needed as well to establish the identity
 \begin{equation}\label{syms}
  g_V(\;A,\;\alpha\,\diamond\,\beta\;)
  \;\;=\;\;
  g_\Z(\;A\,\bullet\,\alpha,\;\beta\;)
  \;\;=\;\;
  g_\Z(\;\alpha,\;A\,\bullet\,\beta\;)
 \end{equation}
 for all $A\,\in\,V$ and $\alpha,\,\beta\,\in\,\Z$, which tells us in
 particular the symmetric spinor multiplication $\diamond$ is completely
 determined by the Clifford multiplication $\bullet$ and thus by the spin
 geometry of $V$. In consequence the spin group $\Spin(\,V,\,-g_V\,)$ acts
 by automorphisms on $\A$ via
 $$
  F\,\star\,(\;a\;\oplus\;\underline a\,\Gamma\;\oplus\;A\;\oplus\;\alpha\;)
  \;\;:=\;\;
  a\;\oplus\;\underline a\,\Gamma\;\oplus\;\Ad_FA\;\oplus\;
  (\,F\,\bullet\,\alpha\,)
 $$
 for all $F\,\in\,\Spin(\,V,-g_V\,)\,\subseteq\,\Cl(\,V,\,-g_V\,)$,
 where $\Ad_F\,\in\,\SO(\,V,\,-g_V\,)\,\stackrel!=\,\SO(\,V,\,+g_V\,)$ denotes
 the restriction of the conjugation by $F$ to the subspace $V\,\subseteq\,
 \Cl(\,V,\,-g_V\,)$. The full automorphism group of $\A$ however is much
 larger due to the presence of the derivations
 \begin{equation}\label{sd}
  \vartheta_\xi(\;a\,\oplus\,\underline a\,\Gamma\,\oplus\,A\,\oplus\;\alpha\;)
  \;\;:=\;\;
  0\;\oplus\;(\,-\,g_\Z(\xi,\alpha)\,)\,\Gamma\;\oplus\;(\,\xi\diamond\alpha\,)
  \;\oplus\;(\,-\,A\bullet\xi\,+\,3\,\underline a\,\xi\,)
 \end{equation}
 for the multiplication (\ref{amult}) parametrized by an arbitrary spinor
 $\xi\,\in\,\Z$. The full automorphism group of the Albert algebra $\A$ turns
 out to be generated by $\Spin(\,V,\,-g_V\,)$ and the exponentials of the
 derivations (\ref{sd}), in particular the automorphism group of $\A$ equals
 the compact simple Lie group $\F_4\,:=\,\Aut\,\A$ of dimension $52$
 \cite{baez} with $\mathbb{Z}_2$--graded Lie algebra
 \begin{equation}\label{f4dec}
  \mathfrak{f}_4
  \;\;:=\;\;
  \mathfrak{aut}\;\A
  \;\;=\;\;
  \so(\,V,\,g_V\,)\;\oplus\;
  \{\;\;\vartheta_\xi\;\;|\;\;\xi\,\in\,\Z\;\;\}
 \end{equation}
 in the sense that the commutator $[\,\vartheta_\xi,\,\vartheta_{\hat\xi}\,]
 \,\in\,\so(\,V,\,g_V\,)$ of two of the derivations (\ref{sd}) lies in
 $\so(\,V,\,g_V\,)$. The symmetric space corresponding to this
 $\mathbb{Z}_2$--graded Lie algebra is the Cayley projective plane
 $\O P^2$, which can be defined for example as the orbit of the element
 $\Gamma\,\in\,\A$ defined in equation (\ref{2el}) under the automorphism
 group $\F_4$ of the Albert algebra $\A$.
 
 \pfill
 For the purpose of this article the most important conclusion of the
 preceeding discussion is that the $26$--dimensional subspace $R\,:=\,
 \Im\,\A\,\subseteq\,\A$ of trace free albertions represents the unique
 non--trivial critical irreducible representation of the simple Lie group
 $\F_4\,=\,\Aut\,\A$ found in Section \ref{critreps}, in fact all other
 non--trivial irreducible representations of $\F_4$ have strictly larger
 dimensions. In particular $R$ decomposes under the stabilizer subgroup
 $\Spin(\,V,\,-g_V\,)\,\cong\,\Spin(\,9\,)$ of the base point
 $\Gamma\,\in\,\Im\,\A$ of the Cayley plane $\O P^2$ into:
 \begin{equation}\label{rdec}
  R
  \;\;:=\;\;
  \mathrm{Im}\;\A
  \;\;=\;\;
  \R\;\oplus\;V\;\oplus\;\Z\ .
 \end{equation}
 According to decomposition (\ref{f4dec}) the isotropy representation of
 the Cayley projective plane $\O P^2$ equals the spinor representation $\Z$
 parametrizing the derivations (\ref{sd}) of the Albert algebra, we may
 thus assume that the Riemannian metric on $\O P^2$ is induced by $g_\Z$.
 Second symmetric and exterior powers of spinor representations $\Z^*\,\cong\,
 \Z$ decompose into exterior powers of the vector representation $V$, in the
 case at hand this decomposition reads:
 \begin{equation}\label{s2dec}
  \S^2\Z^*
  \;\;=\;\;
  \R\;\oplus\;V\;\oplus\;\L^4V
  \qquad\qquad
  \L^2\Z^*
  \;\;=\;\;
  \L^2V\;\oplus\;\L^3V\ .
 \end{equation}
 According to equations (\ref{rdec}) and (\ref{s2dec}) and the Lemma of Schur
 both homomorphism spaces
 $$
  \Hom_{\Spin(V,-g_V)}(\;R,\;\S^2_\circ\Z^*\;)
  \qquad\qquad
  \Hom_{\Spin(V,-g_V)}(\;R,\;\Z\;)
 $$
 are one--dimensional, in particular the Frobenius reciprocal extension
 $F^\ext$ of every non--zero
 $$
  F
  \;\;\in\;\;
  \Hom_{\Spin(V,-g_V)}(\;R,\;\S^2_\circ\Z^*\;)
  \;\;\cong\;\;
  \Hom_{\F_4}(\;R,\;\Gamma(\,\S^2_\circ T^*\O P^2\,)\;)
 $$
 identifies $R$ with the only eigenspace of the Lichnerowicz Laplacian
 $\Delta_L$ on the trace free symmetric $2$--tensors on $\O P^2$ with
 eigenvalue below $2\,\frac\scal{16}$. In light of Corollary \ref{lsrs}
 we need to verify that the prototypical divergence operator $D^*_R$
 associated to $R$ is injective
 \begin{equation}\label{prop2}
  D^*_R:\;\;
  \Hom_{\Spin(V,-g_V)}(\,R,\,\S^2_\circ\Z^*\,)
  \;\longrightarrow\;
  \Hom_{\Spin(V,-g_V)}(\,R,\,\Z\,)
 \end{equation}
 in order to prove that the Cayley projective plane $\O P^2$ is stable. A
 suitable non--zero element $F\,\in\,\Hom_{\Spin(V,-g_V)}(\,R,\,\S^2_\circ
 \Z^*\,)$ is provided by the symmetric spinor multiplication $\diamond$
 $$
  F(\;\underline a\;\oplus\;A\;\oplus\;\alpha\;)
  \;\;:=\;\;
  \frac12\,\sum_{\mu\,\nu}
  g_V(\;A,\;\xi_\mu\,\diamond\,\xi_\nu\;)\;d\xi_\mu\,\cdot\,d\xi_\nu\ ,
 $$
 where $\{\,\xi_\mu\,\}$ is some basis of the isotropy representation $\Z$
 and $\{\,d\xi_\mu\,\}$ denotes the dual basis of $\Z^*$. Using the musical
 isomorphism $\#:\,\Z^*\longrightarrow\Z$ with respect to $g_\Z$ we calculate
 $$
  \tr_{g_\Z}\,F(\;\underline a\;\oplus\;A\;\oplus\;\alpha\;)
  \;\;=\;\;
  g_V\Big(\;A,\;\sum_\mu d\xi^\#_\mu\,\diamond\,\xi_\mu\;\Big)
 $$
 and conclude that the composition $\tr_g\,\circ\,F$ is effectively a
 $\Spin(\,V,\,-g_V\,)$--invariant and thus vanishing linear functional
 on the irreducible representation $V$. In consequence the linear map
 $F:\,R\longrightarrow\S^2_\circ\Z^*$ takes values in the trace free
 symmetric $2$--forms as claimed and:
 \begin{equation}\label{zdia}
  \sum_\mu d\xi^\#_\mu\,\diamond\,\xi_\mu
  \;\;=\;\;
  0\ .
 \end{equation}
 Before proceeding to calculate $D^*_RF$ we want to point out the following
 identity in $\alpha,\,\beta\,\in\,\Z$
 $$
  2\,(\,\beta\,\diamond\,\alpha\,)\,\bullet\,\beta
  \;+\;(\,\beta\,\diamond\,\beta\,)\,\bullet\,\alpha
  \;\;=\;\;
  2\,g_\Z(\,\beta,\,\alpha\,)\,\beta\;+\;g_\Z(\,\beta,\,\beta\,)\,\alpha\ ,
 $$
 which is a partial polarization of the identity $(\alpha\diamond\alpha)
 \,\bullet\,\alpha\,=\,g_\Z(\alpha,\alpha)\,\alpha$ mentioned above.
 Specifically we replace $\alpha$ in equation (\ref{2id}) by the
 expression $\beta+t\alpha$ and take the derivative of the resulting identity
 in $t\,=\,0$ using trilinearity and the symmetry of $\diamond$ and $g_\Z$.
 Tracing the partially polarized identity over $\beta$ with a view on
 equation (\ref{zdia}) we conclude:
 \begin{equation}\label{combs}
  2\,\sum_\mu(\,d\xi^\#_\mu\,\diamond\,\alpha\,)\,\bullet\,\xi_\mu
  \;\;=\;\;
  2\,\sum_\mu g_\Z(\,d\xi^\#_\mu,\,\alpha\,)\,\xi_\mu
  \;+\;\sum_\mu g_\Z(\,d\xi^\#_\mu,\,\xi_\mu\,)\,\alpha
  \;\;=\;\;
  18\,\alpha\ .
 \end{equation}
 In light of the decomposition (\ref{f4dec}) the subspace $\Z\,\subseteq
 \,\mathfrak{f}_4$ of the Lie algebra of $\F_4\,=\,\Aut\,\A$ acts on the
 Albert algebra $\A$ and its invariant subspace $R\,\subseteq\,\A$ by the
 derivations (\ref{sd}), to wit $\xi\,\star\,=\,\vartheta_\xi$ for all
 $\xi\,\in\,\Z$. In consequence the general description (\ref{proto}) of
 the prototypical divergence operator associated to the finite dimensional
 representation $R$ becomes
 \begin{eqnarray*}
  (\,D^*_RF\,)(\;\underline a\;\oplus\;A\;\oplus\;\alpha\;)
  &=&
  \hphantom{\frac12}\,\sum_\lambda
  \Big(\;d\xi^\#_\lambda\,\ins\,F(\;\vartheta_{\xi_\lambda}
  (\,\underline a\,\oplus\,A\,\oplus\,\alpha\,)\;)\;\Big)^\#
  \\
  &=&
  \hphantom{\frac12}\,\sum_\lambda
  \Big(\;d\xi^\#_\lambda\,\ins\,F(\;(\,\ldots\,)\;\oplus\;(\,\xi_\lambda
  \,\diamond\,\alpha\,)\;\oplus\;(\,\ldots\,)\;)\;\Big)^\#
  \\
  &=&
  \frac12\,\sum_{\lambda\,\mu\,\nu}
  g_V(\,\xi_\lambda\,\diamond\,\alpha,\,\xi_\mu\,\diamond\,\xi_\nu\,)\;
  d\xi_\lambda\,\lrcorner\,(\,d\xi^\#_\mu\,\cdot\,d\xi^\#_\nu\,)
  \\
  &=&
  \hphantom{\frac12}\,\sum_{\mu\,\nu}
  g_\Z(\,(\,d\xi^\#_\mu\,\diamond\,\alpha\,)\,\bullet\,\xi_\mu,
  \,\xi_\nu\,)\;d\xi^\#_\nu
  \;\;=\;\;
  9\,\alpha\ ,
 \end{eqnarray*}
 where we have used $g_V(\,A,\,\xi_\mu\diamond\xi_\nu\,)\,=\,g_\Z(\,A\bullet
 \xi_\mu,\,\xi_\nu\,)$ and the trace (\ref{combs}) in the last line. In
 particular the prototypical divergence operator (\ref{prop2}) is injective
 and the Cayley projective plane $\O P^2$ is stable in the sense of Koiso
 as stipulated in Theorem \ref{mt}.
\section{Einstein Deformations of the Grassmannians}
\label{eingr}
 In analogy to our discussion of the stability of the Cayley projective plane
 $\O P^2$ in Section \ref{geocay} we analyze the stability of the family
 of quaternionic Grassmannians $\Gr_r\H^{r+s}$ of $r$--dimensional quaternionic
 subspaces in $\H^{r+s}$ in this section for all parameters $r,\,s\,\geq\,1$.
 Needless to say this family includes the quaternionic projective spaces
 $\H P^s\,=\,\Gr_1\H^{s+1}$ and $S^4\,=\,\H P^1$, whose stability had been
 settled by Koiso except for the case $s\,=\,2$. Calculating the prototypical
 divergence operator we will show that the quaternionic Grassmannians
 $\Gr_r\H^{r+s}$ are stable, if and only if they are quaternionic projective
 spaces in the sense $r\,=\,1$ or $s\,=\,1$.

 \pfill
 For the time being we interprete the unitary symplectic groups $\Sp(n)$
 as the groups of unitary $n\times n$--matrices with coefficients in the
 quaternions $\H$, which act naturally on the spaces $\H^n$ of column vectors
 considered as right vector spaces over $\H\,=\,\mathrm{Mat}_{1\times 1}\H$
 by right matrix multiplication. For all $r,\,s\,\geq\,1$ the induced action of
 $\Sp(r+s)$ on the Grassmannian
 $$
  \Gr_r\H^{r+s}
  \;\;=\;\;
  \Sp(\,r\,+\,s\,)/_{\displaystyle\Sp(\,r\,)\,\times\,\Sp(\,s\,)}
 $$
 of $r$--dimensional subspaces of $\H^{r+s}$ is transitive as well with
 stabilizer in the base point $\H^r\,\subseteq\,\H^{r+s}$ given by the
 diagonal subgroup $\Sp(r)\,\times\,\Sp(s)\,\subseteq\,\Sp(r+s)$. Under
 restriction to the stabilizer of the base point the defining representation
 $V\,:=\,\H^{r+s}$ decomposes into the direct sum $V\,=\,H\,\oplus\,E$ of the
 defining representations $H\,:=\,\H^r$ and $E\,:=\,\H^s$ of $\Sp(r)$ and
 $\Sp(s)$. In general we will consider $V$ as a complex vector space via
 the obvious inclusion $\C\,\subseteq\,\H$ endowed with the
 $\Sp(r+s)$--equivariant, conjugate linear map $C:\,V\longrightarrow V$ of
 right multiplication with $j$ satisfying $C^2\,=\,-\id_V$; analogous remarks
 apply to both $H$ and $E$. The adjoint representation of the symplectic Lie
 groups equals the second symmetric power of its defining representation, in
 turn the decomposition $V\,=\,H\,\oplus\,E$ implies
 \begin{eqnarray*}
  \sp(\,r+s\,)\,\otimes_\R\C
  \;\;=\;\;
  \S^2V
  &=&
  \S^2H\,\oplus\,(\,H\,\otimes\,E\,)\,\oplus\,\S^2E
  \\
  &=&
  \sp(\,r\,)\,\otimes_\R\C\,\oplus\,(\,H\,\otimes\,E\,)
  \,\oplus\,\sp(\,s\,)\,\otimes_\R\C
 \end{eqnarray*}
 so that $H\otimes E\,\cong\,T_{\H^r}M\,\otimes_\R\C$ equals the
 complexified isotropy representation of the irreducible symmetric space
 $M\,:=\,\Gr_r\H^{r+s}$. Its trace free second symmetric power
 decomposes into:
 \begin{equation}\label{grdec}
  \S^2_\circ(\,H\otimes E\,)
  \;\;=\;\;
  (\,\S^2H\otimes\S^2E\,)
  \;\oplus\;(\,\L^2_\circ H\otimes\L^2_\circ E\,)
  \;\oplus\;(\,\L^2_\circ H\,)\;\oplus\;(\,\L^2_\circ E\,)\ .
 \end{equation}
 En nuce the all important difference between the projective spaces with
 $\min\{\,r,\,s\,\}\,=\,1$ and the general Grassmannians with $r,\,s\,\geq\,2$
 derives from this decomposition, after all $\L^2_\circ H\,=\,\{\,0\,\}$ or
 $\L^2_\circ E\,=\,\{\,0\,\}$ vanish for $r\,=\,1$ or $s\,=\,1$ respectively.

 According to the classification of critical representations of the
 symplectic Lie algebras $\sp(2n,\C)\,=\,\sp(n)\otimes_\R\C$ in
 Section \ref{critreps} the only critical irreducible representations
 of $\sp(r+s)$ besides the trivial and adjoint representations are
 $V$, $\L^2_\circ V$ together with $\L^3_\circ V$ for $r+s\,=\,3$.
 In light of the decomposition (\ref{grdec}) the Frobenius reciprocity
 (\ref{frob}) ensures that neither $V$ nor $\L^3_\circ V$ occurs in the
 trace free symmetric $2$--tensors $\Gamma(\,\S^2_\circ T^*M\,)$, and
 so we are left with
 \begin{equation}\label{r2dec}
  R\otimes_\R\C
  \;\;:=\;\;
  \L^2_\circ V
  \;\;=\;\;
  \C\;\oplus\;(\,\L^2_\circ H\,)\;\oplus\;(\,H\,\otimes\,E\,)
  \;\oplus\;(\,\L^2_\circ E\,)
 \end{equation}
 as the only candidate, where $R\,\subseteq\,\L^2_\circ V$ equals the
 subspace of real elements with respect to the equivariant real structure.
 Using Schur's Lemma we conclude from (\ref{grdec}) and (\ref{r2dec}) that
 \begin{eqnarray*}
  \Hom_{\Sp(r+s)}(\;R,\;
  \hbox to74pt{\hfill$\Gamma(\,\S^2_\circ T^*M\,)$\hfill}\;)\,\otimes_\R\C
  &\cong&
  \Hom_{\Sp(r)\,\times\,\Sp(s)}(\;\L^2_\circ V,\;
  \hbox to74pt{\hfill$\S^2_\circ(\,H\,\otimes\,E\,)$\hfill}\;)
  \\[3pt]
  \Hom_{\Sp(r+s)}(\;R,\;
  \hbox to74pt{\hfill$\Gamma(\,TM\,)$\hfill}\;)\,\otimes_\R\C
  &\cong&
  \Hom_{\Sp(r)\,\times\,\Sp(s)}(\;\L^2_\circ V,\;
  \hbox to74pt{\hfill$H\,\otimes\,E$\hfill}\;)
 \end{eqnarray*}
 have dimension $2$ and $1$ respectively unless $M$ is a quaternionic
 projective space. In consequence the prototypical divergence operator
 associated to the critical representation $\L^2_\circ V$
 \begin{equation}\label{grproto}
  D^*_{\L^2_\circ V}:\;\;\Hom_{\Sp(r)\times\Sp(s)}
  (\;\L^2_\circ V,\;\S^2_\circ(\,H\,\otimes\,E\,)\;)
  \;\longrightarrow\;\Hom_{\Sp(r)\times\Sp(s)}(\;\L^2_\circ V,\;H\otimes E\;)
 \end{equation}
 can not be injective in general, because its domain has dimension $2$ and
 is range dimension $1$. According to Corollary \ref{lsrs} the Einstein
 metrics on the quaternionic Grassmannians $\Gr_r\H^{r+s}$ with parameters
 $r,\,s\,\geq\,2$ are thus unstable in the sense of Koiso.

 \pfill
 In order to study the prototypical divergence operator $D^*_{\L^2_\circ V}$
 in somewhat more detail we recall that $H$ and $E$ come along with complex
 bilinear, invariant symplectic forms $\sigma_H$ and $\sigma_E$ respectively,
 which give rise to the musical isomorphisms $\b:\,H\longrightarrow H^*,\,
 h\longmapsto\sigma_H(\,h,\,\cdot\,)$ and $\#\,:=\,\b^{-1}$ as well as the
 analogous musical isomorphisms $\b$ and $\#$ for $E$. Choosing pairs of
 dual bases $\{\,h_\alpha\,\}$, $\{\,dh_\alpha\,\}$ and $\{\,e_\mu\,\}$,
 $\{\,de_\mu\,\}$ for $H$ and $E$ we define the two linear maps
 \begin{eqnarray*}
  F^H(\;a\,\wedge\,\hat a\;)
  &=&
  \sum_\mu
  (a^\b\otimes de_\mu)\;\cdot\;(\hat a^\b\otimes e^\b_\mu)
  \;-\;\frac{\sigma(\,a,\,\hat a\,)}{2\,r}\;\sum_{\alpha\mu}
  (dh_\alpha\otimes de_\mu)\,\cdot\,(h^\b_\alpha\otimes e^\b_\mu)
  \\
  F^E(\;f\,\wedge\,\hat f\;)
  &=&
  \sum_\alpha
  (dh_\alpha\otimes f^\b)\,\cdot\,(h^\b_\alpha\otimes\hat f^\b)
  \;-\;\frac{\sigma(\,f,\hat f\,)}{2\,s}\;\sum_{\alpha\mu}
  (dh_\alpha\otimes de_\mu)\,\cdot\,(h^\b_\alpha\otimes e^\b_\mu)
 \end{eqnarray*}
 from $\L^2H$ and $\L^2E$ respectively to $\S^2(H\otimes E)^*$. Both
 maps are evidently equivariant under $\Sp(r)\,\times\,\Sp(s)$ and
 send the bivectors $\sum_\alpha dh^\#_\alpha\wedge h_\alpha$ and
 $\sum_\mu de^\#_\mu\wedge e_\mu$ to zero. In essence they thus
 factorize through linear maps $\L^2_\circ H\longrightarrow\S^2_\circ
 (H\otimes E)^*$ and $\L^2_\circ E\longrightarrow\S^2_\circ(H\otimes E)^*$,
 which we may precompose with the corresponding projections in (\ref{r2dec})
 to linear maps
 $$
  F^H:\;\;\L^2_\circ V\;\longrightarrow\;\S^2_\circ(\,H\,\otimes\,E\,)^*
  \qquad\qquad
  F^E:\;\;\L^2_\circ V\;\longrightarrow\;\S^2_\circ(\,H\,\otimes\,E\,)^*\ ,
 $$
 which generate the domain $\Hom_{\Sp(r)\times\Sp(s)}(\,\L^2_\circ F,
 \,\S^2_\circ(H\otimes E)^*\,)$ of $D^*_{\L^2_\circ V}$ due to the
 decomposition (\ref{grdec}). Needless to say the special case of
 quaternionic projective spaces $\min\{\,r,\,s\,\}$ is reflected by
 $F^H\,=\,0$ or $F^E\,=\,0$ for $r\,=\,1$ or $s\,=\,1$. Without loss
 of generality we may assume that the symplectic forms $\sigma_H$ and
 $\sigma_E$ are normalized to make their tensor product agree with the
 Riemannian metric $g\,=\,\sigma_H\otimes\sigma_E$ in the base point
 $\H^r\,\in\,M$ so that
 $$
  \sum_\lambda E_\lambda\,\otimes\,E_\lambda
  \;\;=\;\;
  \sum_{\alpha\mu}(\,dh^\#_\alpha\,\otimes\,de^\#_\mu\,)\,\otimes\,
  (\,h_\alpha\,\otimes\,e_\mu\,)
 $$
 holds true for every orthonormal basis $\{\,E_\lambda\,\}$ of
 $T_{\H^r}M\,\otimes_\R\C\,\cong\,H\otimes E$. Replacing the sum
 over the orthonormal basis $\{\,E_\lambda\,\}$ in formula (\ref{proto})
 by the sum on the right hand side of this identity we calculate
 in a first step
 \begin{eqnarray*}
  (\,D^*_{\L^2_\circ V}F^H\,)(\,a\otimes f\,)^\b
  &=&
  +\;\sum_{\alpha\mu}(\,dh^\#_\alpha\otimes de^\#_\mu\,)\,\ins\,F^H\bigg(\;
  (\,h_\alpha\otimes e_\mu\,)\;\star_{\L^2_\circ V}\;(\,a\wedge f\,)\;\bigg)
  \\
  &=&
  +\;\sum_{\alpha\mu}(\,dh^\#_\alpha\otimes de^\#_\mu\,)
  \,\ins\,F^H\bigg(\,\sigma_H(\,h_\alpha,a\,)\,e_\mu\wedge f
  \,+\,\sigma_E(\,e_\mu,f\,)\,a\wedge h_\alpha\,\bigg)
  \\
  &=&
  -\;\sum_\alpha(\,dh^\#_\alpha\otimes f\,)\,\ins\,F^H(\;a\wedge h_\alpha\;)
  \;\;=\;\;
  \sum_\alpha(\,h_\alpha\otimes f\,)\,\ins\,F^H(\;a\wedge dh^\#_\alpha\;)
 \end{eqnarray*}
 due to the antisymmetry $\sum dh^\#_\alpha\otimes h_\alpha\,=\,-\sum
 h_\alpha\otimes dh^\#_\alpha$ and the definition of the representation
 $\star_{\L^2_\circ V}$ of the Lie algebra $\sp(r+s)\otimes_\R\C\,=\,\S^2V$.
 Using this intermediate result we conclude:
 \begin{eqnarray*}
  (\,D^*_{\L^2_\circ V}F^H\,)(\,a\otimes f\,)^\b
  &=&
  +\;\sum_{\alpha\mu}(\,h_\alpha\otimes f\,)\;\ins\;\bigg(\;
  (\,a^\b\otimes de_\mu\,)\,\cdot\,(\,dh_\alpha\otimes e^\b_\mu\,)\;\bigg)
  \\
  &&
  +\;\sum_{\alpha\beta\mu}\frac{\sigma_H(\,dh^\#_\alpha,a\,)}{2\,r}
  \;(\,h_\alpha\otimes f\,)\;\ins\;\bigg(\;(\,dh_\beta\otimes de_\mu\,)
  \,\cdot\,(\,h^\b_\beta\otimes e^\b_\mu\,)\;\bigg)
  \\
  &=&
  +\;\Big(\;1\;-\;2\,r\;+\;\frac1r\;)\,(\,a^\b\otimes f^\b\,)
  \;\;=\;\;
  -\;\frac{(r-1)(2r+1)}r\;(\,a^\b\,\otimes\,f^\b\,)\ .
 \end{eqnarray*}
 The completely analogous calculation for the linear map $F^E$ results in
 $$
  (\,D_{\L^2_\circ V}^*F^E\,)(\;a\,\otimes\,f\;)
  \;\;=\;\;
  +\;\frac{(s-1)(2s+1)}s\;(\,a\,\otimes\,f\,)\ ,
 $$
 where the sign change is effectively caused by the reversed order of the
 product $e_\mu\wedge f$ compared to $a\wedge h_\alpha$ in the formula for
 $(\,h_\alpha\otimes e_\mu\,)\,\star_{\L^2_\circ V}\,(\,a\wedge f\,)$. In
 consequence the prototypical divergence operator (\ref{grproto}) always
 has maximal rank for all $r,\,s\,\geq\,1$, in particular $D^*_{\L^2_\circ V}$
 is injective for $r\,=\,1$ or $s\,=\,1$ proving the stability of the Einstein
 metric on all quaternionic projective spaces including $S^4\,=\,\H P^1$ and
 $\H P^2$ in the sense of Koiso.
\section*{Acknowledgements}
 The first author acknowledges the support received by the Special Priority
 Program SPP 2026 {\em Geometry at Infinity} funded by the Deutsche
 Forschungsgemeinschaft DFG. Likewise the second author expresses his
 gratitude for the funding received as a SNI member from the Consejo
 Nacional de Ciencia y Tecnolog\'ia CONACyT.
\end{document}